\documentclass[11pt]{amsart}

\usepackage{amsmath,amssymb,amsthm,amsrefs,a4wide}
\usepackage{latexsym}
\usepackage{indentfirst}
\usepackage{graphicx}
\usepackage{placeins}
\usepackage{booktabs}
\usepackage{algorithm}
\usepackage{algorithmic}
\usepackage{multirow}
\usepackage{color}

\theoremstyle{plain}

\theoremstyle{definition}

\newtheorem{example}{Example}

\theoremstyle{remark}
\newtheorem{remark}{Remark}

\newcommand{\diag}{\mathrm{diag}}

\newcommand{\bbm}{\begin{bmatrix}}
\newcommand{\ebm}{\end{bmatrix}}
\newcommand{\R}{\mathbb{R}}
\newcommand{\C}{\mathbb{C}}
\newcommand{\D}{\mathbb{D}}

\newcommand{\Z}{\mathbb{Z}}
\renewcommand{\S}{\mathbb{S}}

\newcommand{\ut}{\tilde{u}}
\newcommand{\xt}{\tilde{x}}
\newcommand{\wt}{\tilde{w}}
\renewcommand{\tt}{\tilde{t}}
\newcommand{\pt}{\tilde{p}}

\newcommand{\Zt}{\tilde{Z}}

\newcommand{\Gh}{\hat{G}}

\newcommand{\bg}{\mathbf{g}}
\newcommand{\bv}{\mathbf{v}}
\newcommand{\bu}{\mathbf{u}}
\newcommand{\bp}{\mathbf{p}}

\newcommand{\but}{\tilde{\bu}}

\newcommand{\bpt}{\tilde{\bp}}

\newcommand{\bgh}{\hat{\bg}}

\begin{document}

\title[Eigenmatrix]{Eigenmatrix for unstructured sparse recovery}

\author[]{Lexing Ying} \address[Lexing Ying]{Department of Mathematics, Stanford University,
  Stanford, CA 94305} \email{lexing@stanford.edu}

\thanks{This work is partially supported by NSF grants DMS-2011699 and DMS-2208163. The author thanks Laurent Demanet and Lin Lin for helpful discussions and the anonymous referee for constructive comments.}



\begin{abstract}
  This note considers the unstructured sparse recovery problems in a general form. Examples include rational approximation, spectral function estimation, Fourier inversion, Laplace inversion, and sparse deconvolution. The main challenges are the noise in the sample values and the unstructured nature of the sample locations. This note proposes the eigenmatrix, a data-driven construction with desired approximate eigenvalues and eigenvectors. The eigenmatrix offers a new way for these sparse recovery problems. Numerical results are provided to demonstrate the efficiency of the proposed method.
\end{abstract}

\maketitle

%----------------------------------------------------------
\section{Introduction}\label{sec:intro}

%---------
This note considers the unstructured sparse recovery problems of a general form. Let $X$ be the parameter space, typically a subset of $\R$ or $\C$, and $S$ be the sampling space. $G(s,x)$ is the kernel function for $s\in S$ and $x\in X$, and is assumed to be analytic in $x$. Suppose that
\[
f(x) = \sum_{k=1}^{n_x} w_k \delta(x-x_k)
\]
is the unknown sparse signal, where $n_x$ is the number of spikes, $\{x_k\}$ are the spike locations, and $\{w_k\}$ are the spike weights. The observable of the problem is
\[
u(s) :=\int_X G(s,x) f(x) dx = \sum_{k=1}^{n_x} G(s,x_k) w_k
\]
for $s\in S$.

Let $\{s_j\}$ be a set of $n_s$ unstructured sample locations in $S$ and $u_j:=u(s_j)$ be the exact values. Suppose that we are only given the noisy observations $\ut_j:=u_j (1+\sigma Z_j)$, where $Z_j$ are independently identically distributed (i.i.d) random variables with zero mean and unit variance, and $\sigma$ is the noise magnitude. The task is to recover the spike locations $\{x_k\}$ and weights $\{w_k\}$.

Quite a few sparse recovery problems can be cast into this general form. Below is a partial list.
\begin{itemize}
\item Rational approximation. $G(s,x) = \frac{1}{s-x}$, $X$ is typically a set in $\C$, and $\{s_j\}$ are locations separated from $X$. Two common cases of $X$ are the unit disk and the half-plane.

\item Spectral function estimation of many-body quantum systems. $G(s,x)= \frac{1}{s-x}$, $X$ is a real interval $[-b,b]$ of the complex plane, and $\{s_j\}$ is the Matsubara grid on the imaginary axis.

\item Fourier inversion. For example $G(s,x) = \exp(\pi i sx)$, $X$ is the interval $[-1,1]$, and $\{s_j\}$ is a set of real numbers.

\item Laplace inversion. $G(s,x) = x \exp(-sx)$, $X$ is an interval $[c_1,c_2]$ of the positive real axis, and $\{s_j\}$ is a set of positive real numbers.

\item Sparse deconvolution. $G(s,x)$ is a translational invariant kernel, such as $G(s,x) = \frac{1}{1+ \gamma (s-x)^2}$, $X$ is a real interval, and $\{s_j\}$ is a set of real numbers.
\end{itemize}

%Sparse recovery problems are known to be difficult.
The primary challenges of the current setup come from two sources. First, the kernel $G(s,x)$ can be quite general. Second, the sample locations $\{s_j\}$ are unstructured, which excludes many existing algorithms that exploit special structures. Third, the sample values $\{\ut_j\}$ are noisy, which raises stability issues when the recovery problem is quite ill-posed.

%---------
\subsection{Contribution}
This note introduces the {\em eigenmatrix} for these unstructured sparse recovery problems. By defining the vector-valued function $\bg(x)=[G(s_j,x)]_{1\le j\le n_s}$ for $x\in X$, we introduce the eigenmatrix as an $n_s\times n_s$ matrix $M$ that satisfies, for any $x\in X$
\[
M \bg(x) \approx x \bg(x).
\]
This is a a data-driven object that depends on $G(\cdot,\cdot)$, $X$, and the sample locations $\{s_j\}$. The main features of the eigenmatrix are
\begin{itemize}
\item It assumes no special structure of the sample locations $\{s_j\}$.
\item It offers a rather unified approach to these sparse recovery problems.
\item As the numerical results suggest, even when the recovery problem is ill-conditioned, the reconstruction can be quite robust with respect to noise.
\end{itemize}

%---------
\subsection{Related work} There has been a long list of works devoted to the sparse recovery problems mentioned above.

Rational approximation has a long history in numerical analysis. Some of the well-known methods are the RKFIT algorithm \cite{berljafa2017rkfit}, barycentric interpolation \cite{berrut2004barycentric}, Pade approximation \cite{gonnet2013robust}, vector fitting \cite{gustavsen1999rational}, and AAA \cite{nakatsukasa2018aaa}.

Spectral function approximation is a key computational task for many-body quantum systems. Well-known methods include Pade approximation \cite{vidberg1977solving,beach2000reliable,schott2016analytic}, maximum entropy methods \cite{jarrell1996bayesian,beach2004identifying,levy2017implementation,kraberger2017maximum,rumetshofer2019bayesian}, and stochastic analytic continuation \cite{sandvik1998stochastic,vafayi2007analytical,goulko2017numerical,krivenko2019triqs}. Several most recent algorithms are \cite{fei2021nevanlinna,fei2021analytical,ying2022analytic,ying2022pole,huang2023robust}.

Fourier inversion is a vast field with many different problem setups. When $X$ and $S$ are dual discrete grids with $\{s_j\}$ chosen randomly, this is the compressive sensing problem \cite{candes2006stable,donoho2006most,foucart2013invitation}, and there is a vast literature on methods based on the $\ell_1$ convex relaxation. When $X$ is an interval and $\{s_j\}$ are equally spaced grid points, this becomes the line spectrum estimation or superresolution problem \cite{donoho1992superresolution}. Both Prony-type methods \cite{prony1795essai,schmidt1986multiple,roy1989esprit,hua1990matrix} and optimization-based approaches \cite{candes2014towards,demanet2015recoverability,moitra2015super} are well-studied for this field.

Laplace inversion is a longstanding computational problem. Most established algorithms \cite{abate2004multi,weideman1999algorithms,weeks1966numerical,weideman2007parabolic} assume the capability of accessing the sample values at any arbitrary locations. For the case of equally-spaced sample locations, Prony-type methods have been proposed in \cite{beylkin2009nonlinear,potts2013parameter}. The work in \cite{peter2013generalized,stampfer2020generalized} further extends the Prony-type methods to the kernels associated with more general first-order and second-order differential operators.

For sparse deconvolution, when $\{s_j\}$ forms a uniform grid, it is closely related to the superresolution problem. However, when $\{s_j\}$ are unstructured, the literature is surprisingly limited.

The rest of the note is organized as follows. Section \ref{sec:pe} reviews Prony's method and the ESPRIT algorithms for the special case of the exponential kernel with the uniform sampling grid. Section \ref{sec:em} describes the eigenmatrix approach for the general kernels and unstructured grids. Section \ref{sec:num} presents the numerical experiments of the applications mentioned above. Section \ref{sec:disc} concludes with a discussion for future work.

%----------------------------------------------------------
\section{Prony and ESPRIT}\label{sec:pe}

To motivate the eigenmatrix construction, we first briefly review Prony's method and the ESPRIT algorithm. Consider the recovery problem with $X=\S\equiv\{z: |z|=1\}\subset\C$, $S=\R$, $G(s,x) = x^s$ (with the branch cut at $x=-1$), and $s_j=j$ for $j\in\Z$. Here, we make the simplifying assumption that $\{s_j\}$ is the whole integer lattice, though only a finite chunk is required in the actual implementations. Most presentations of Prony's method and the ESPRIT algorithm start with the Hankel matrix. However, our presentation emphasizes the role of the shifting operator in order to motivate the eigenmatrix approach in Section \ref{sec:em}.

Introduce the infinitely-long vector-valued function $\bg(x)=[G(s_j,x)]_{j\in\Z} = [x^j]_{j\in\Z}$. Let $M$ be the shifting operator that moves each entry up by one slot, i.e., $(M \bv)_j = \bv_{j+1}$ for any vector $\bv$.  Then
\[
M \bg(x) = x \bg(x). %\quad \text{or equivalently}\quad (M - e^{ i x}I) \bg(x) = 0. 
\]
Define the vector $\bu=[u_j]_{j\in\Z}$ of the exact observations $u_j=u(s_j) = \sum_k G(s_j,x_k) w_k = \sum_k w_k x_k^j$. Since $\bu=\sum_k \bg(x_k) w_k$, we have
for any $t\ge 0$
\[
M^t \bu = \sum_k \bg(x_k) w_k x_k^t.
\]

For the Prony's method, consider the matrix
\[
\begin{bmatrix}  \bu & M\bu & \ldots & M^{n_x} \bu \end{bmatrix}
=
\bbm \bg(x_1) & \ldots & \bg(x_{n_x}) \ebm
\bbm w_1 & &\\& \ddots & \\& & w_{n_x}\ebm
\bbm
1 &  x_1 & \ldots & x_1^{n_x} \\
\vdots & \vdots & \ddots & \vdots \\
1 &  x_{n_x} & \ldots & x_{n_x}^{n_x}
\ebm
\]
Let $\bp$ be a non-zero vector in its null space, i.e.,
\begin{equation}
  \begin{bmatrix}  \bu & M\bu & \ldots & M^{n_x} \bu \end{bmatrix} \bp =0 \quad\text{with}\quad
  \bp = \begin{bmatrix} p_0 \\ \vdots\\ p_{n_x} \end{bmatrix} \label{eq:solvep}.
\end{equation}
Therefore,
\[
\bbm
1 &  x_1 & \ldots & x_1^{n_x} \\
\vdots & \vdots & \ddots & \vdots \\
1 &  x_{n_x} & \ldots & x_{n_x}^{n_x}
\ebm
\bbm p_0 \\ \vdots \\ p_{n_x} \ebm = 0.
\]
This implies that $\{x_k\}$ are the roots of $p(x)=p_0 + p_1 x + \ldots p_{n_x} x^{n_x}$. Therefore, one can identify $\{x_k\}$ via rootfinding once $\bp$ is computed.

For the ESPRIT algorithm, consider the matrix
\[
\begin{bmatrix}  \bu & M\bu & \ldots & M^{\ell} \bu \end{bmatrix}
=
\bbm \bg(x_1) & \ldots & \bg(x_{n_x}) \ebm
\bbm w_1 & &\\& \ddots & \\& & w_{n_x}\ebm
\bbm
1 &  x_1    & \ldots & x_1^\ell \\
\vdots & \vdots & \ddots & \vdots \\
1 &  x_{n_x} & \ldots & x_{n_x}^\ell
\ebm
\]
with $\ell>n_x$. Let $USV^*$ be the rank-$n_x$ singular value decomposition (SVD) of this matrix. The matrix $V^*$ takes the form
\[
V^* = P
\bbm
1 &  x_1    & \ldots & x_1^\ell \\
\vdots & \vdots & \ddots & \vdots \\
1 &  x_{n_x} & \ldots & x_{n_x}^\ell
\ebm
\]
where $P$ is an unknown non-degenerate $n_x \times n_x$ matrix. Let $Z_0$ and $Z_1$ be the submatrices of $V^*$ by excluding the first column and the last column, respectively, i.e.,
\[
Z_0 =
P
\bbm
1      & \ldots & x_1^{\ell-1} \\
\vdots & \ddots & \vdots \\
1      & \ldots & x_{n_x}^{\ell-1}
\ebm,
\quad
Z_1 =
P
\bbm
x_1    &  \ldots & x_1^\ell \\
\vdots     &  \ddots & \vdots \\
x_{n_x} & \ldots & x_{n_x}^\ell
\ebm.
%% = 
%% P
%% \bbm
%% e^{ i x_1} &  & \\
%% & \ddots & \\
%% & & e^{ i x_{n_x}}
%% \ebm
%% Z_0.
\]
By introducing
\[
Z_1(Z_0)^+ =
P
\bbm
x_1 &  & \\
& \ddots & \\
& & x_{n_x}
\ebm
P^{-1},
\]
one can identify $\{x_k\}$ by computing the eigenvalues of $Z_1 (Z_0)^+$.

For both methods, given $\{x_k\}$, the sample weights $\{w_k\}$ can be computed via, for example, the least square solve.

%% Since $(M - e^{ i x_\ell}I) \bg(x_\ell) = 0$ for each $\ell$ and the operators $M - e^{ i x_\ell}I$ for different $x_\ell$ commute with each other, we have
%% \[
%% \prod_{\ell=1}^{n_x} (M- e^{ i x_\ell} I) \bg(x_k) = 0.
%% \]
%% Define the vector $\bu=[u_j]_{j\in\Z}$ of the exact observations $u_j=u(s_j) = \sum_k G(s_j,x_k) w_k = \sum_k e^{ i j x_k} w_k = \sum_k w_k \bg(x_k)$. Since $\bu$ is a linear combination of $\{\bg(x_k)\}_{1\le k \le n_x}$, we have
%% \begin{equation}
%%   \prod_{\ell=1}^{n_x} (M- e^{ i x_\ell} I) \bu = 0. \label{eq:prodM}
%% \end{equation}
%% Because $\prod_\ell (M- e^{ i x_\ell} I)$ is a matrix polynomial of form $p(M) := p_0 + p_1 M + \cdots + p_{n_x} M^{n_x}$,
%% \eqref{eq:prodM} translates to 
%% \[
%% (p_0 + p_1 M + \cdots + p_{n_x} M^{n_x}) \bu=0.
%% \]
%% This is equivalent to
%% \begin{equation}
%%   \begin{bmatrix}  \bu & M\bu & \ldots & M^{n_x} \bu \end{bmatrix} \bp =0 \quad\text{with}\quad
%%   \bp = \begin{bmatrix} p_0 \\ \vdots\\ p_{n_x} \end{bmatrix} \label{eq:solvep}
%% \end{equation}
%% One can solve for the vector $\bp$ by identifying the null space vector of the matrix in \eqref{eq:solvep}. The roots of the polynomial $p(x)=p_0 + p_1 x + \ldots p_{n_x} x^{n_x}$ are $\{e^{ i x_k}\}$. 

\begin{remark}
  For most problems, the sample values $\{\ut_j\}$ have noise. As a result, the sample locations and weights obtained above are only approximations. Many implementations of the Prony and ESPRIT methods have a postprocessing step, where these approximations are used as the initial guesses of the following optimization problem
  \begin{equation}
    \min_{\xt_k,\wt_k} \sum_j \left|\sum_k \xt_k^j \wt_k - \ut_j\right|^2. \label{eq:optwk}
  \end{equation}
\end{remark}

\begin{remark}
  For an actual problem, the number of spikes $n_x$ is not known a priori. An important question is how to pick the right degree of the polynomial $p(x)$ (for Prony's method) or the rank of the truncated SVD (for the ESPRIT algorithm). The general criteria are that the objective value of \eqref{eq:optwk} (after postprocessing) should be within the noise level, and the degree $d$ should be as small as possible. Commonly used criteria include AIC \cite{akaike1998information} and BIC \cite{schmidt1986multiple}.
\end{remark}

%----------------------------------------------------------
\section{Eigenmatrix}\label{sec:em}

%============
\subsection{Main idea}

The discussion above uses two special features of the problem: (a) the kernel is of the exponential form, and (b) $\{s_j\}$ forms an equally spaced grid. These two features together allow one to write down the shifting operator $M$ explicitly. However, these two features no longer hold for sparse recovery problems with a general kernel $G(\cdot,\cdot)$ or unstructured sample locations $\{s_j\}$.

To address this challenge, we take a data-driven approach.  Let $\bg(x)$ now be the $n_s$-dimensional vector $[G(s_j,x)]_{1\le j \le n_s}$. The main idea is to introduce an {\em eigenmatrix} $M$ of size $n_s\times n_s$ such that for all $x\in X$
\[
M \bg(x) \approx x \bg(x).
\]
The reason why $M$ is called the eigenmatrix is because {\em it is designed to have the desired approximate eigenvalues and eigenvectors}.

Below, we detail how to apply the eigenmatrix idea to complex and real cases. Since the choice of degree (in Prony's algorithm) and the rank of truncated SVD (in the ESPRIT algorithm) remain unchanged, we present the algorithm for fixed $n_x$ in order to simplify the discussion.

\subsection{Complex analytic case}

To simplify the discussion, assume first that $X$ is the unit disc $\D$, and we will comment on the general case at the end. Define for each $x$ the vector $\bg(x):=[G(s_j,x)]_{1\le j\le n_s}$. The first step is to construct $M$ such that $M \bg(x) \approx x \bg(x)$ for $x\in \D$. Numerically, it is more robust to use the normalized vector $\bgh(x) = \bg(x)/\|\bg(x)\|$ since the norm of $\bg(x)$ can vary significantly depending on $x$. The condition then becomes
\[
M \bgh(x) \approx x \bgh(x), \quad x\in \D.
\]
We enforce this condition on a uniform grid $\{a_t\}_{1\le t \le n_a}$ of size $n_a$ on the boundary of the unit disk
\[
M \bgh(a_t) \approx a_t \bgh(a_t).
\]
Define an $n_s\times n_a$ matrix $\Gh = [\bgh(a_t)]_{1\le t\le n_a}$ with $\bgh(a_t)$ as columns and also an $n_a\times n_a$ diagonal matrix $\Lambda = \diag(a_t)$. The above condition can be written in a matrix form as
\begin{equation}
  M \Gh \approx \Gh \Lambda. \label{eq:MG}
\end{equation}
\begin{remark}
  The main guideline for the choice of $n_a$ is that the columns of $\Gh$ are numerically linearly independent. To see why this is essential, let us consider the extreme case of a kernel $G(s,x)$ constant in $x$. Here, the columns are linearly dependent, and there is no way to recover individual $\{x_k\}$. In practice, different sampling locations $\{s_j\}$ leads to different choices of $n_a$. In practice, $n_a$ is chosen such that the condition number of $\Gh$ is bounded below $10^7$.
\end{remark}

When the columns of $\Gh$ are numerically linearly independent, \eqref{eq:MG} suggests the following choice of the eigenmatrix
\begin{equation}
  M := \Gh \Lambda \Gh^\dagger, \label{eq:M}
\end{equation}
where the pseudoinverse $\Gh^\dagger$ is computed by thresholding the singular values of $\Gh$. Since $M$ is to be applied repetitively as in \eqref{eq:solvep}, the threshold is set so that the norm of $M$ is bounded by a small constant such as 3.

\begin{remark}
  A key question is why enforcing the condition on the uniform grid $\{a_t\}$ is enough. The following calculation shows why.  $M \bg(a_t) \approx a_t \bg(a_t)$ implies $M \bg(x) \approx x \bg(x)$ for all $x\in \D$:
  \begin{align*}
    M \bg(x) & = \frac{1}{2\pi i} \int_{\partial\D} \frac{M \bg(z)}{z-x} dz \approx \frac{1}{2\pi i}\sum_t \frac{M \bg(a_t)}{a_t-x} ( i a_t \frac{2\pi}{n_a})\\
    & \approx \frac{1}{2\pi i}\sum_t \frac{a_t \bg(a_t)}{a_t-x} (i a_t \frac{2\pi}{n_a}) \approx \frac{1}{2\pi i} \int_{\partial\D} \frac{z \bg(z)}{z-x} dz
    = x \bg(x),
  \end{align*}
  where the first and third approximations use the exponential convergence of the trapezoidal rule for analytic functions $\bg(x)$ and $x\bg(x)$, and the second approximation directly comes from $M \bg(a_t) \approx a_t \bg(a_t)$. The equalities are applications of the Cauchy integral theorem.
\end{remark}

\begin{remark}
  The second question is, what if $X$ is not $\D$? For a general connected domain $X$ with smooth boundary, let $\phi(t): \D \rightarrow X$ be the one-to-one map between $\D$ and $X$ from the Riemann mapping theorem. We then consider the new kernel $G(s,t) = G(s, \phi(t))$ between $S$ and $\D$ and use the above algorithm to recover the locations $\{\tt_k\}$ in $\D$. Once $\{\tt_k\}$ are available, we set $\xt_k =\phi(\tt_k)$.
\end{remark}

\begin{figure}[h!]
  \centering
  \includegraphics[scale=0.25]{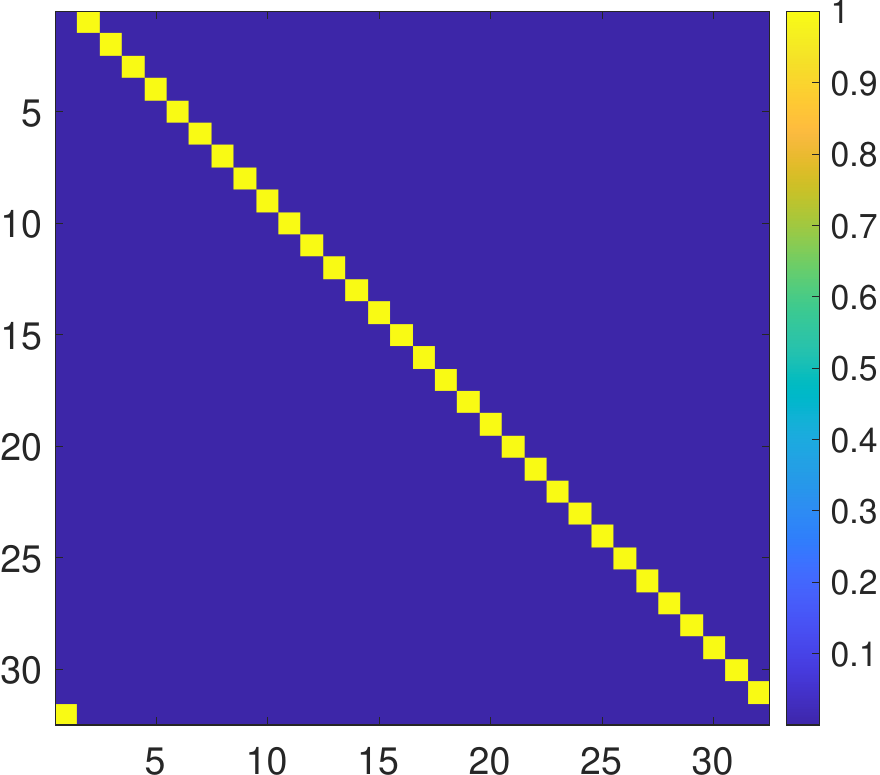}
  \includegraphics[scale=0.25]{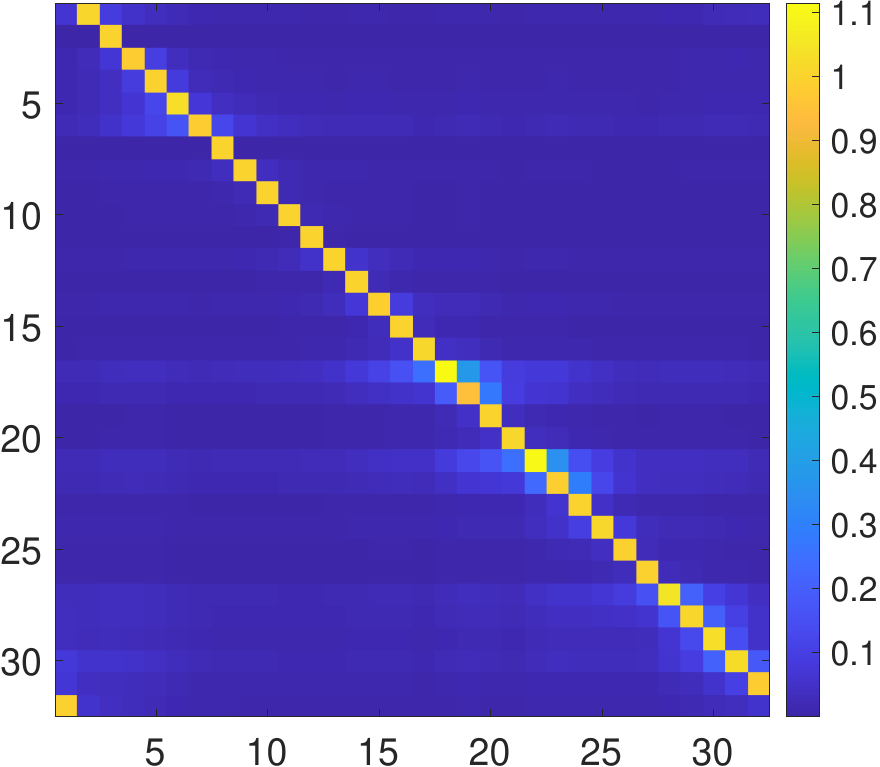}
  \includegraphics[scale=0.25]{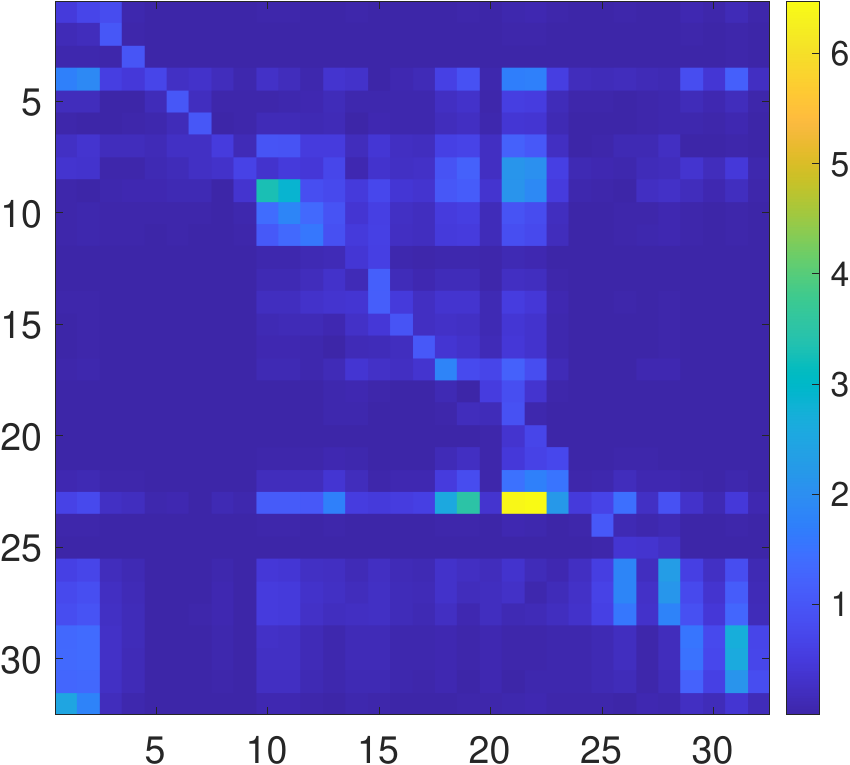}
  \caption{The eigenmatrix for the example in Section \ref{sec:pe}. $X=\S\equiv\{z: |z|=1\}\subset\C$, $S=\R$, $G(s,x) = x^s$ (with the branch cut at $x=-1$). $n_s=32$.  Left: $M$ when $s_j=j$ for $0\le j <n_s$. Middle: $M$ when $s_j$ is a random perturbation of the integer lattice. Right: $M$ when $\{s_j\}$ are chosen uniformly in $[0,n_s]$.}
  \label{fig:cmp}
\end{figure}

To provide an idea about the eigenmatrix, we apply the above construction to the problem given in Section \ref{sec:pe} with $n_s=32$. Figure \ref{fig:cmp} plots the eigenmatrix for three cases: (a) $\{s_j=j\}$ is a integer lattice (the classical case); (b) $\{s_j\}$ is a random perturbation of the integer lattice; (c) $\{s_j\}$ are chosen uniformly in $[0,n_s]$. In the integer lattice case, the eigenmatrix approach reproduces the shifting matrix. In the perturbative case, the eigenmatrix is quite close to the shifting matrix, as expected. In the fully random case, the eigenmatrix is quite different but still keeps the overall trend.

\subsection{Real analytic case}
To simplify the discussion, assume that $X$ is the interval $[-1,1]$, and we will comment on the general case later. Let us define for each $x$ the vector $\bg(x) = [G(s_j,x)]_{1\le j\le n_s}$. The first step is to construct $M$ such that $M \bg(x) \approx x \bg(x)$ for $x\in [-1,1]$. Numerically, it is again more robust to use the normalized vector $\bgh(x) = \bg(x)/\|\bg(x)\|$ and consider the modified condition
\[
M \bgh(x) \approx x \bgh(x), \quad x\in [-1,1].
\]
We enforce this condition on a Chebyshev grid $\{a_t\}_{1\le t \le n_a}$ of size $n_a$ on the interval $[-1,1]$:
\[
M \bgh(a_t) \approx a_t \bgh(a_t).
\]
Introduce the $n_s\times n_s$ matrix $\Gh = [\bgh(a_t)]_{1\le t\le n_a}$ with columns $\bgh(a_t)$ as well as the $n_a\times n_a$ diagonal matrix $\Lambda = \diag(a_t)$. The condition now reads
\[
M \Gh \approx \Gh \Lambda.
\]
When the columns of $\Gh$ are numerically linearly independent, this again suggests the following choice of the eigenmatrix for the real analytic case
\[
M := \Gh \Lambda \Gh^\dagger,
\]
where the pseudoinverse $\Gh^\dagger$ is computed by thresholding the singular values of $\Gh$.

\begin{remark}
  We claim that, for real analytic kernels $G(s,x)$, enforcing the condition at the Chebyshev grid $\{a_t\}$ is sufficient. To see this,
  \[
  M \bg(x) \approx M\left(\sum_k c_t(x) \bg(a_t)\right) = \sum_t c_t(x) M \bg(a_t) \approx \sum_t c_t(x) (a_t \bg(a_t)) \approx x \bg(x),
  \]
  where $c_t(x)$ is the Chebyshev quadrature for $x$ associated with grid $\{a_t\}$. Here, the first and third approximations use the convergence property of the Chebyshev quadrature for analytic functions $\bg(x)$ and $x\bg(x)$, and the second approximation directly comes from $M \bg(a_t) \approx a_t \bg(a_t)$.
\end{remark}

\begin{remark}
  The next question is, what if $X$ is not the interval $[-1,1]$? For a general interval or analytic segment $X$, let $\phi: [-1,1] \rightarrow X$ be a smooth one-to-one map between $[-1,1]$ and $X$. By considering the kernel $G(s,t)= G(s, \phi(t))$ instead and applying the above algorithm, we can recover the locations $\{\tt_k\}$ in $[-1,1]$. Finally, set $\xt_k=\phi(\tt_k)$.
\end{remark}

%============
\subsection{Putting together}

%% With $M$ available, we define the vector $\but:=[\ut_j]_{1\le j \le n_s}$, form the system
%% \[
%% \begin{bmatrix}  \but & M\but & \ldots & M^d \but \end{bmatrix} \bpt = 0 \quad\text{with}\quad
%% \bpt = \begin{bmatrix} \pt_0 \\ \vdots\\ \pt_d \end{bmatrix},
%% \]
%% and solve for $\bpt$ as the null space vector. The roots of the polynomial $\pt(x) := \pt_0 + \pt_1 x + \cdots + \pt_d x^d$ are the estimators $\{\xt_k\}$ for the spike locations and a least square solve gives $\{\wt_k\}$ as the estimators for the spike weights.

With the eigenmatrix $M$ available, the rest is similar to Prony's method and the ESPRIT algorithm.
Define the vector $\but:=[\ut_j]_{1\le j\le n_s}$ from the noisy sample values.

For Prony's method, consider
\[
\bbm  \but & M\but & \ldots & M^{n_x} \but \ebm
\approx
\bbm \bg(x_1) & \ldots & \bg(x_{n_x}) \ebm
\bbm w_1 & &\\& \ddots & \\& & w_{n_x}\ebm
\bbm
1 &  x_1 & \ldots & x_1^{n_x} \\
\vdots & \vdots & \ddots & \vdots \\
1 &  x_{n_x} & \ldots & x_{n_x}^{n_x}
\ebm.
\]
Let $\bpt$ be a non-zero vector in its null-space 
\[
\bbm  \but & M\but & \ldots & M^{n_x} \but \ebm \bpt =0 \quad\text{with}\quad \bpt =
\bbm \pt_0 \\ \vdots\\ \pt_{n_x} \ebm
\]
Therefore,
\[
\bbm
1 &  x_1 & \ldots & x_1^{n_x} \\
\vdots & \vdots & \ddots & \vdots \\
1 &  x_{n_x} & \ldots & x_{n_x}^{n_x}
\ebm
\bbm \pt_0 \\ \vdots \\ \pt_{n_x} \ebm \approx 0.
\]
The roots of the polynomial $\pt(x) := \pt_0 + \pt_1 x + \cdots + \pt_{n_x} x^{n_x}$ provide the estimators $\{\xt_k\}$ for $\{x_k\}$. 

For the ESPRIT algorithm, consider the matrix
\[
\bbm  \but & M\but & \ldots & M^\ell \but \ebm
\approx
\bbm \bg(x_1) & \ldots & \bg(x_{n_x}) \ebm
\bbm w_1 & &\\& \ddots & \\& & w_{n_x}\ebm
\bbm
1 &  x_1 & \ldots & x_1^\ell \\
\vdots & \vdots & \ddots & \vdots \\
1 &  x_{n_x} & \ldots & x_{n_x}^\ell
\ebm
\]
with $\ell>n_x$. Let $\tilde{U} \tilde{S}\tilde{V}^*$ be the rank-$n_x$ truncated SVD of this matrix. The matrix $\tilde{V}^*$ satisfies
\[
\tilde{V}^* \approx P
\bbm
1 &  x_1 & \ldots & x_1^\ell \\
\vdots & \vdots & \ddots & \vdots \\
1 &  x_{n_x} & \ldots & x_{n_x}^\ell
\ebm
\]
where $P$ is an unknown non-degenerate $n_x \times n_x$ matrix. Let $\Zt_0$ and $\Zt_1$ be the submatrices of $\Zt$ by excluding the first column and the last column, respectively, i.e.,
\[
\Zt_0 \approx P
\bbm
1      & \ldots & x_1^{\ell-1} \\
\vdots & \ddots & \vdots \\
1      & \ldots & x_{n_x}^{\ell-1}
\ebm,
\quad
\Zt_1 \approx P
\bbm
x_1 & \ldots & x_1^\ell \\
\vdots & \ddots & \vdots \\
x_{n_x} & \ldots & x_{n_x}^\ell
\ebm.
%% \approx
%% P
%% \bbm
%% x_1 &  & \\
%% & \ddots & \\
%% & & x_{n_x}
%% \ebm
%% \Zt_0
\]
By introducing
\[
\Zt_1 (\Zt_0)^+ \approx
P
\bbm
x_1 &  & \\
& \ddots & \\
& & x_{n_x}
\ebm
P^{-1},
\]
one can get estimates $\{\xt_k\}$ for $\{x_k\}$ by computing the eigenvalues of $\Zt_1 (\Zt_0)^+$.

With $\{\xt_k\}$ available, the least square solve
\[
\min_{\wt_k} \sum_j \left|\sum_k G(s_j,\xt_k) \wt_k - \ut_j\right|^2 
\]
gives the estimators $\{\wt_k\}$ for $\{w_k\}$ for both methods.

%% Since $(M-x_k I) \bg(x_k)\approx 0$,
%% \[
%% \prod_\ell (M-x_\ell I) \bg(x_k)\approx 0 \quad\text{and}\quad  \prod_\ell(M-x_\ell I) \but \approx 0.
%% \]
%% Therefore, we can set up the following system:
%% \[
%% \begin{bmatrix}  \but & M\but & \ldots & M^d \but \end{bmatrix} \bpt = 0 \quad\text{with}\quad
%% \bpt = \begin{bmatrix} \pt_0 \\ \vdots\\ \pt_d \end{bmatrix},
%% \]
%% and solve for $\bpt$ as the null space vector. 

%----------------------------------------------------------
\section{Numerical results}\label{sec:num}

This section applies the eigenmatrix approach to the unstructured sparse recovery problems mentioned in Section \ref{sec:intro}. In all examples, the spike weights $\{w_k\}$ are set to be $1$ and the noises $\{Z_j\}$ are Gaussian.

Once the eigenmatrix $M$ is constructed, the reported numerical results are obtained using ESPRIT. The results of Prony's method are similar but slightly less robust. In each plot, blue, green, and red stand for the exact solution, the result before postprocessing, and the one after postprocessing.

%-------
\begin{example}[Rational approximation] The problem setup is 
  \begin{itemize}
  \item $G(s,x) = \frac{1}{s-x}$.
  \item $X=\D$.
  \item $\{s_j\}$ are random points outside the disk, each with a modulus between $1.2$ and $2.2$. $n_s=40$.
  \end{itemize}

  \begin{figure}[h!]
    \centering
    \includegraphics[scale=0.25]{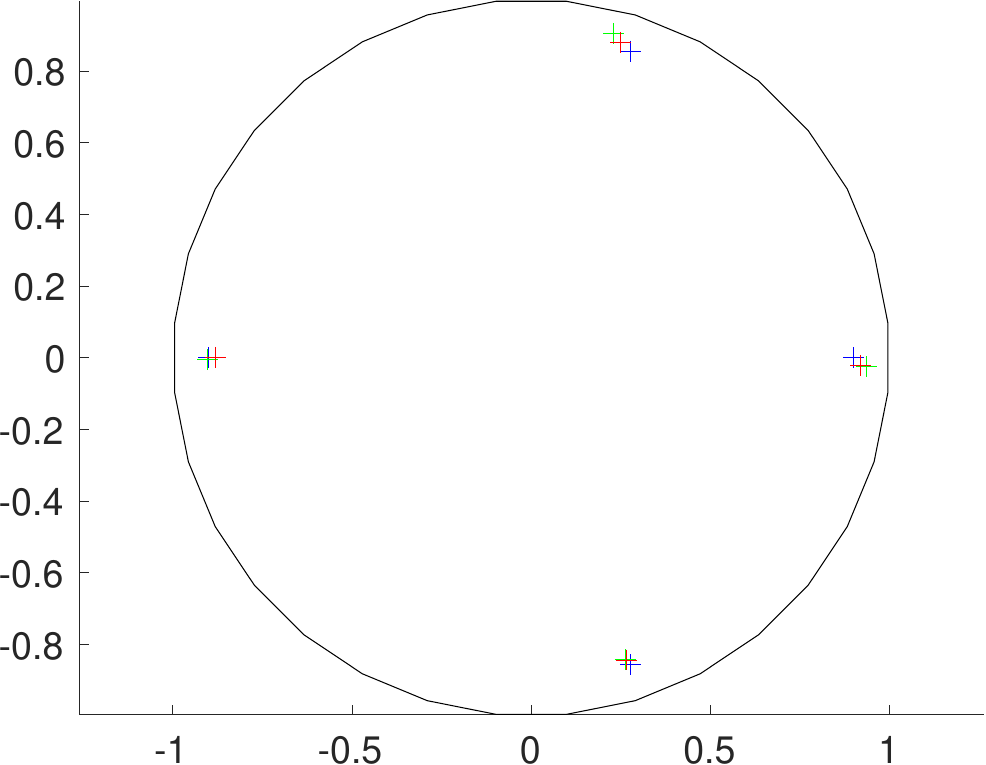}
    \includegraphics[scale=0.25]{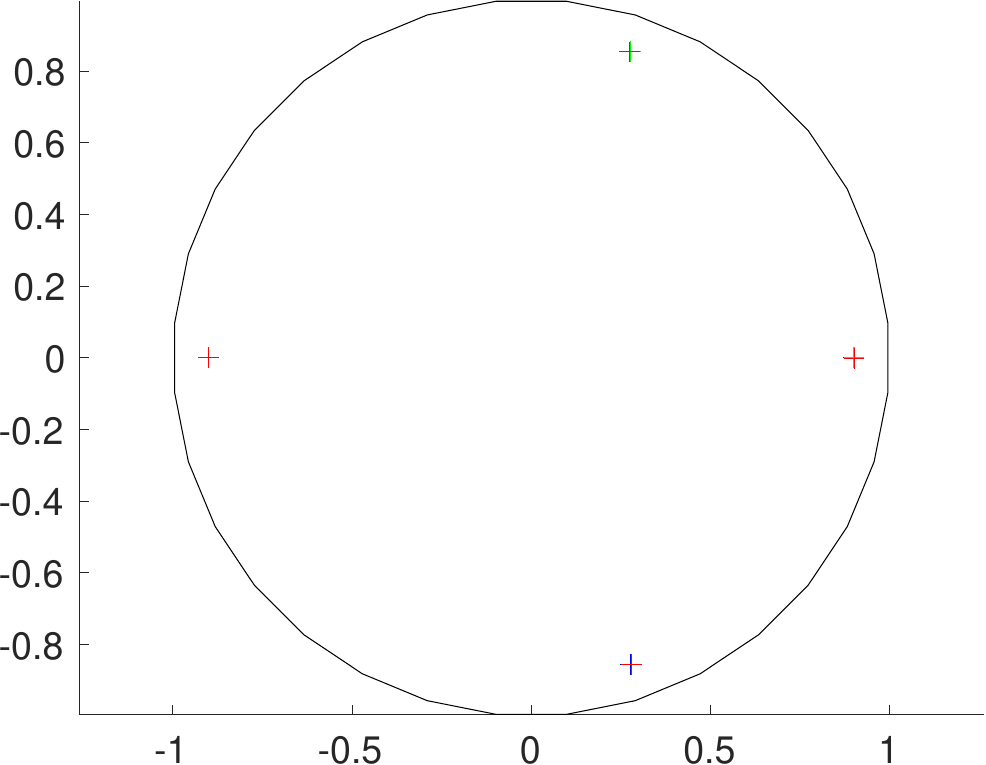}
    \includegraphics[scale=0.25]{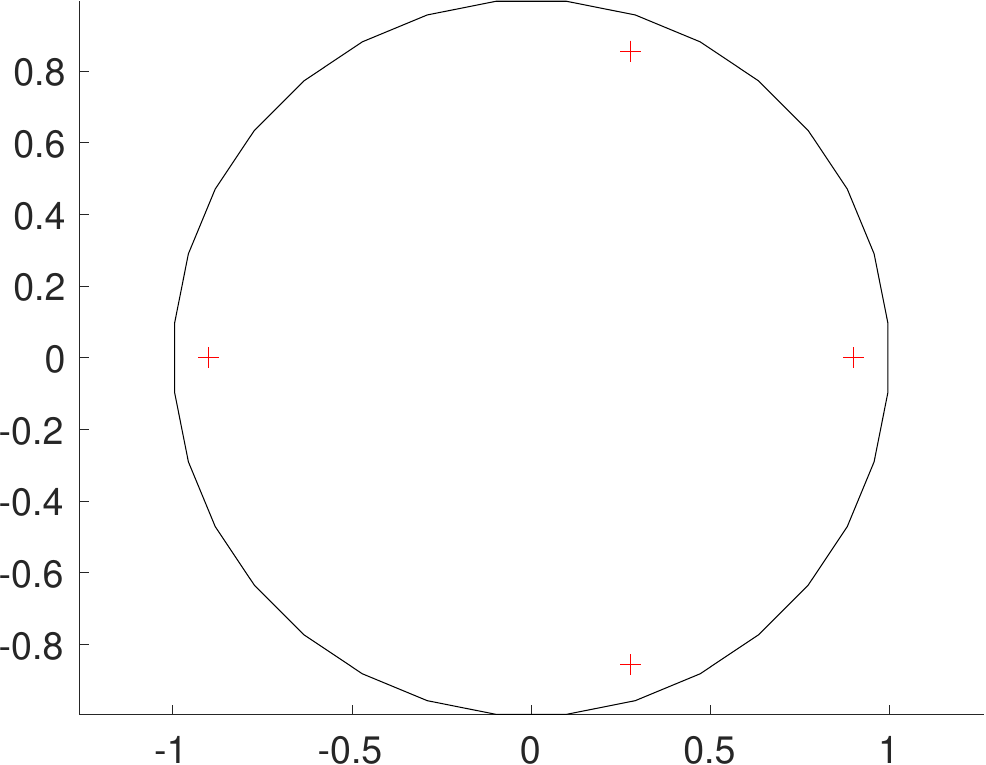}\\
    \includegraphics[scale=0.25]{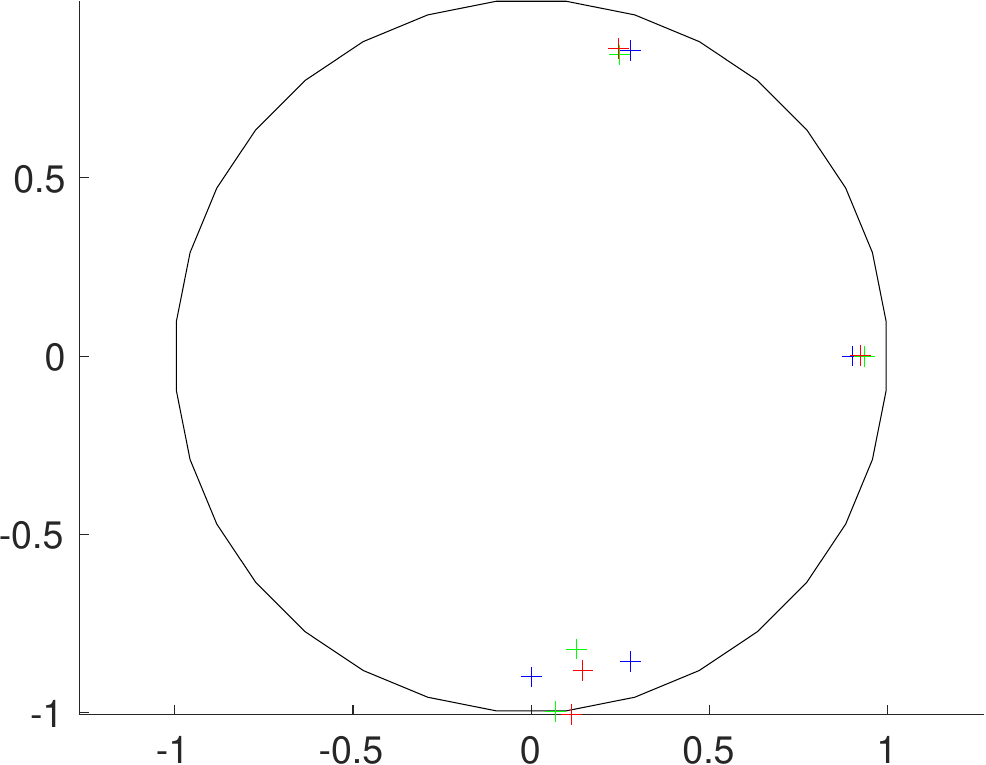}
    \includegraphics[scale=0.25]{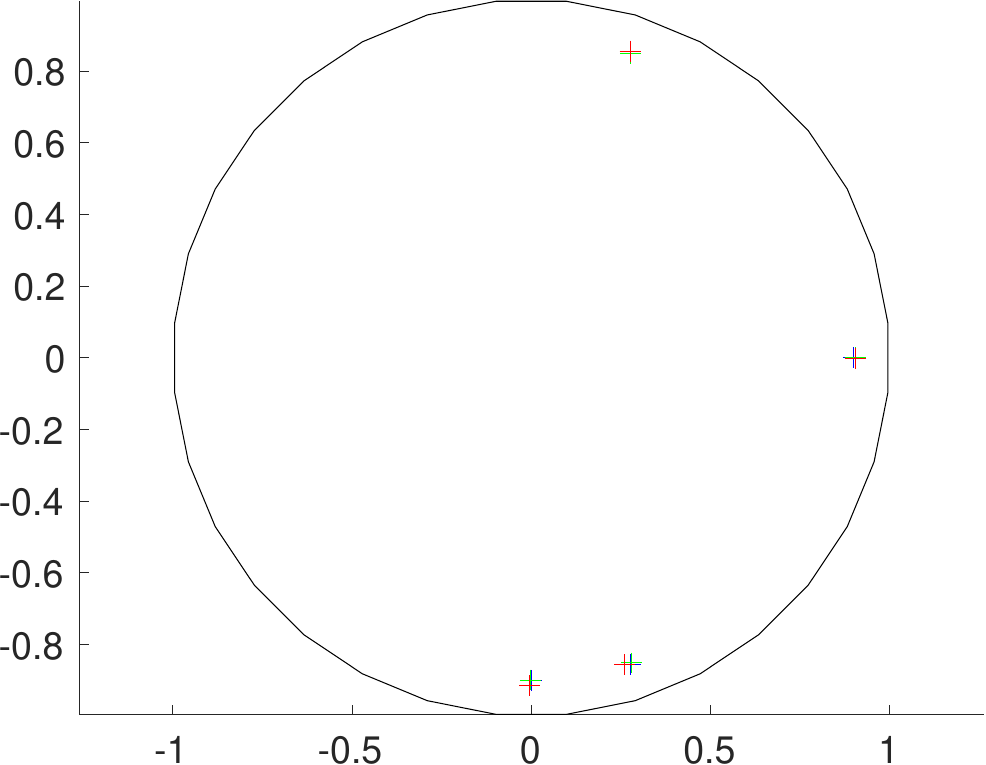}
    \includegraphics[scale=0.25]{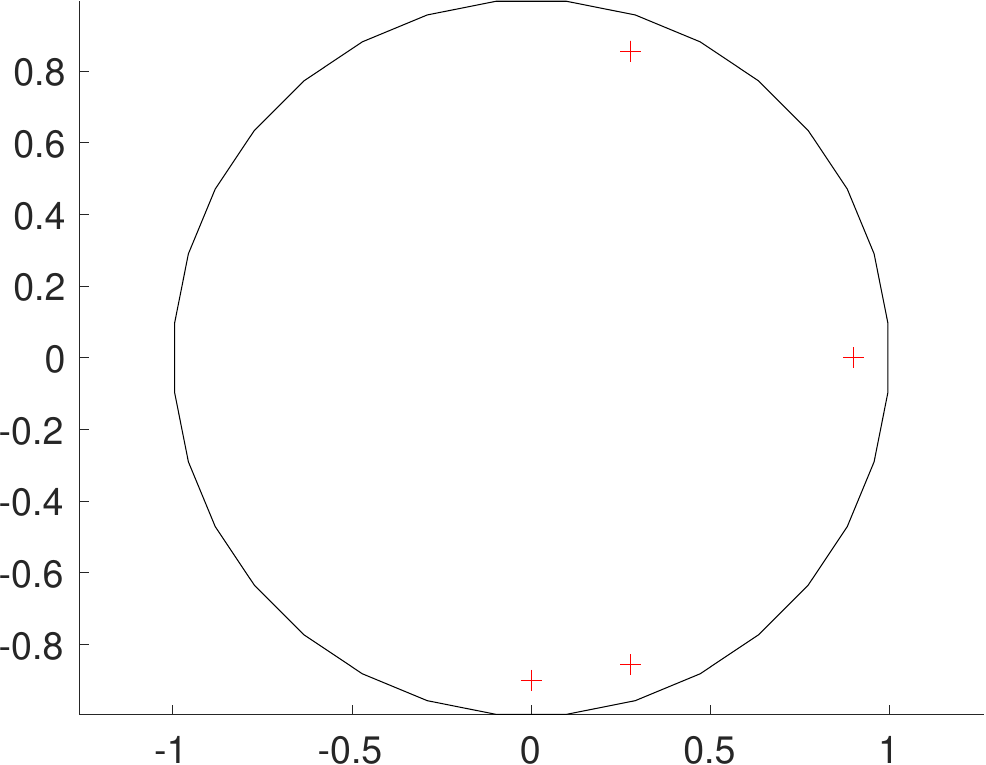} 
    \caption{Rational approximation. $G(s,x) = \frac{1}{s-x}$.  $X=\D$. $\{s_j\}$ are random points outside the unit disk, each with a modulus between $1.2$ and $2.2$. $n_s=40$. Columns: $\sigma$ equals to $10^{-2}$, $10^{-3}$, and $10^{-4}$. Rows: the easy test (well-separated) and the hard test (with two nearby spikes).}
    \label{fig:exR}
  \end{figure}
  
  Figure \ref{fig:exR} summarizes the experimental results. $n_a=32$. The three columns correspond to noise levels $\sigma$ equal to $10^{-2}$, $10^{-3}$, and $10^{-4}$.  Two tests are performed. In the first one (top row), $\{x_k\}$ are well-separated from each other. The plots show accurate recovery of the spike locations from all $\sigma$ values. In the second one (bottom row), two of the spike locations are close to each other: the reconstruction at $\sigma=10^{-2}$ shows a noticeable error, while the results for $\sigma=10^{-3}$ and $\sigma=10^{-4}$ are accurate. The plots also suggest that the eigenmatrix approach results in fairly accurate (green) initial guesses for the postprocessing step.
\end{example}

%-------
\begin{example}[Spectral function approximation] The problem setup is
  \begin{itemize}
  \item $G(s,x) = \frac{1}{s-x}$.
  \item $X=[-1,1]$.
  \item $\{s_j\}$ is the Matsubara grid from $-\frac{(2N-1)\pi}{\beta}i$ to $\frac{(2N-1)\pi}{\beta}i$ with $\beta=100$ and $N=128$. Hence, $n_s=256$.
  \end{itemize}

  \begin{figure}[h!]
  \centering
  \includegraphics[scale=0.25]{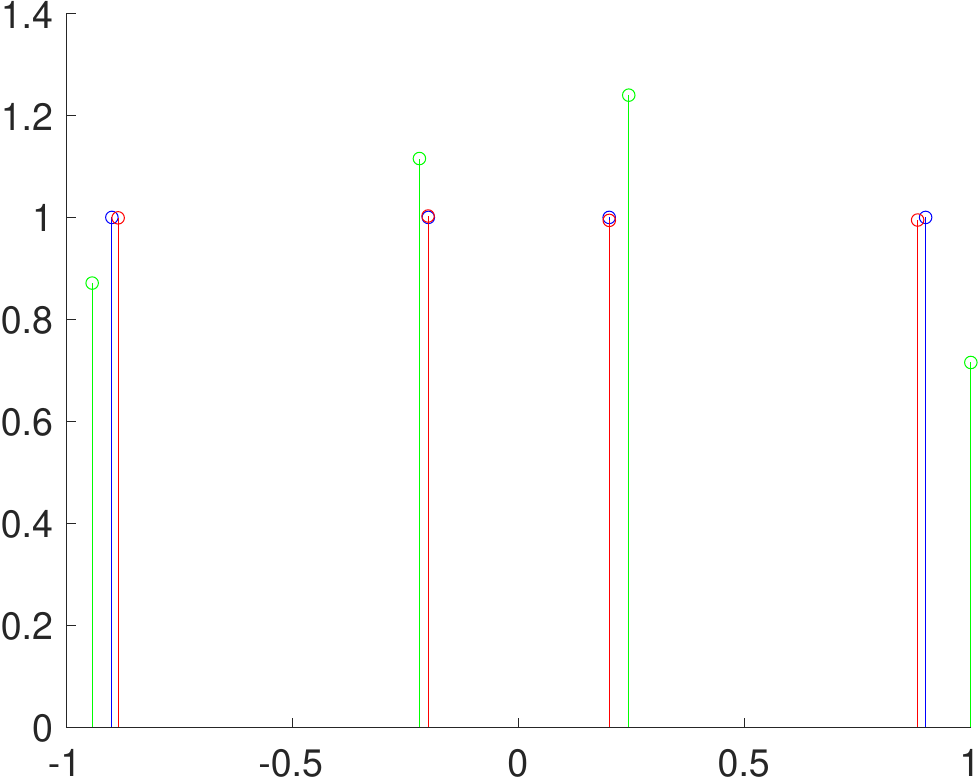}
  \includegraphics[scale=0.25]{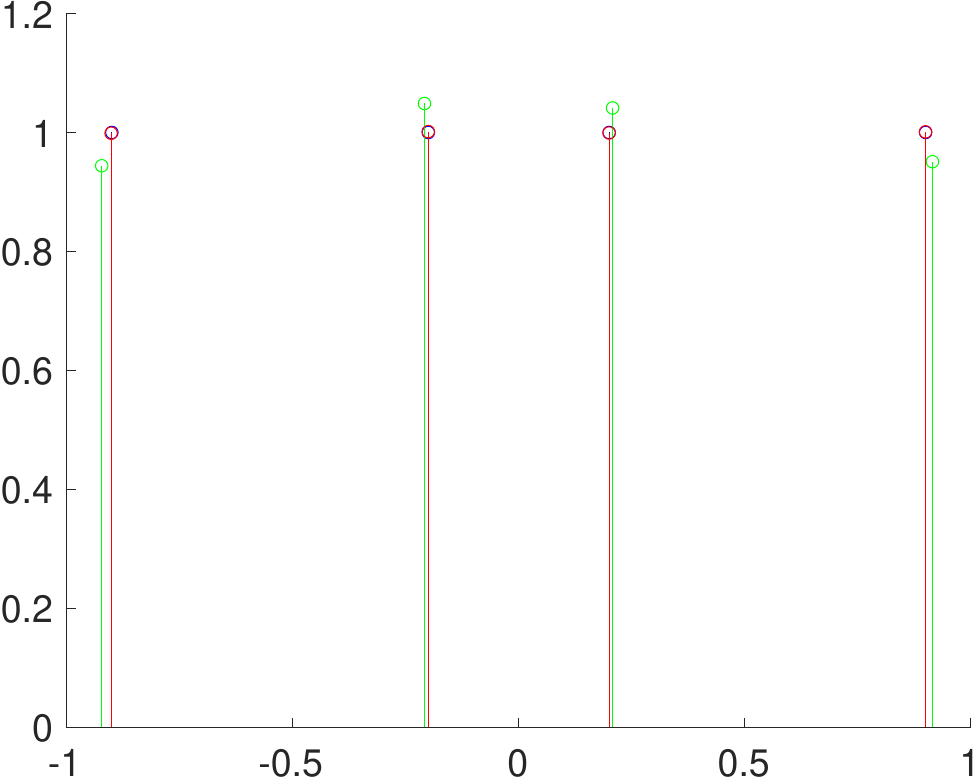}
  \includegraphics[scale=0.25]{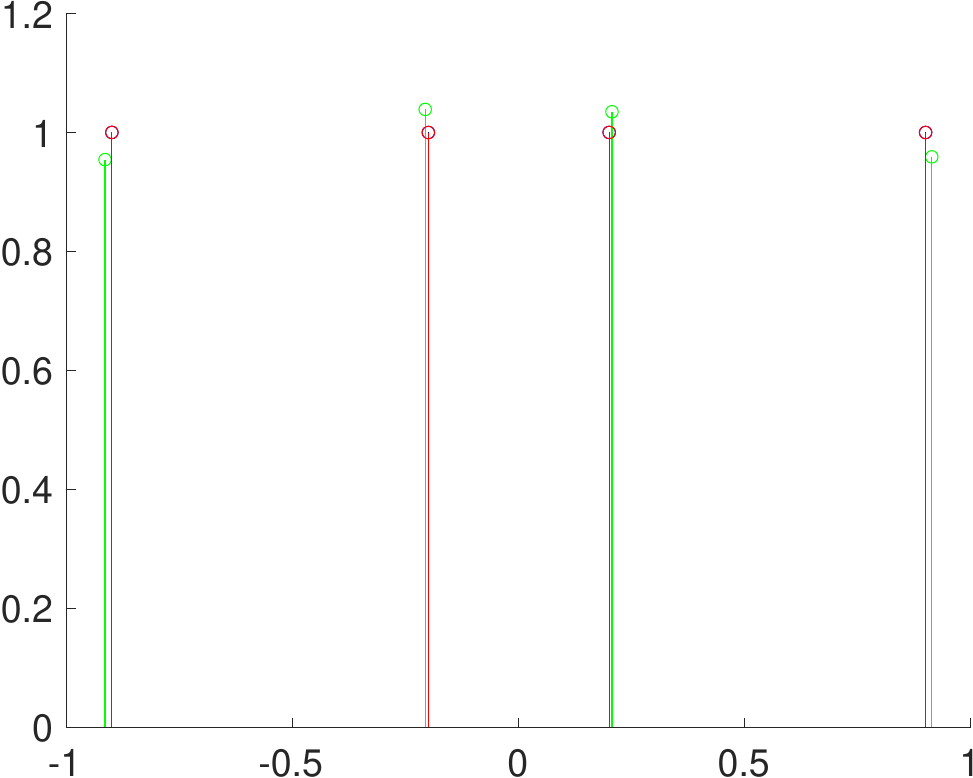}\\
  \includegraphics[scale=0.25]{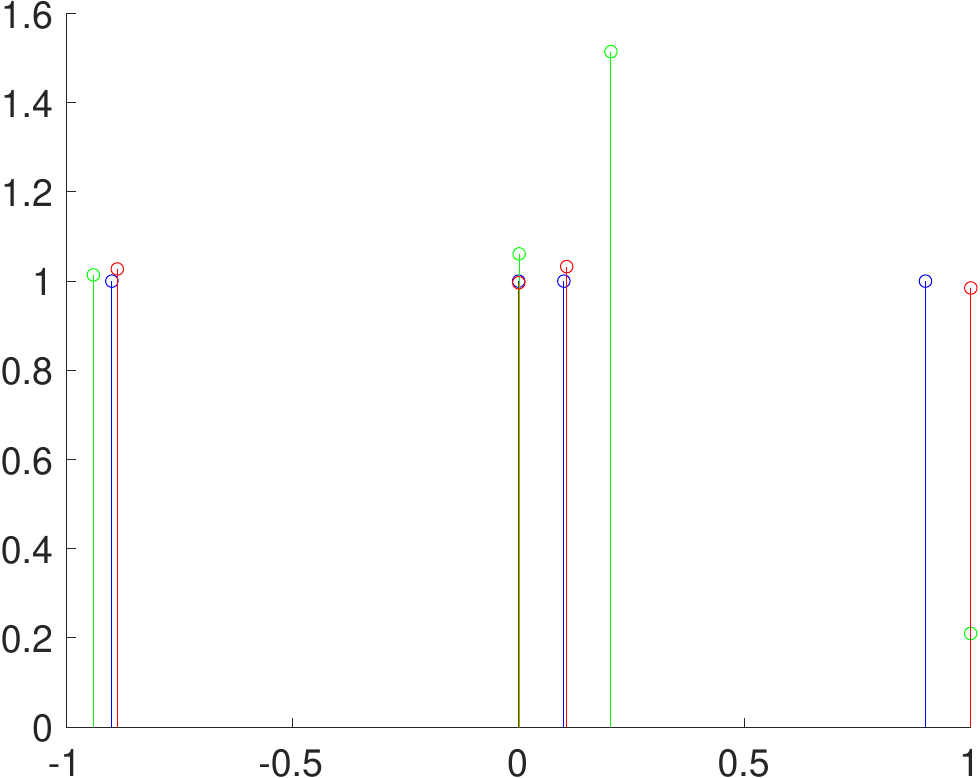}
  \includegraphics[scale=0.25]{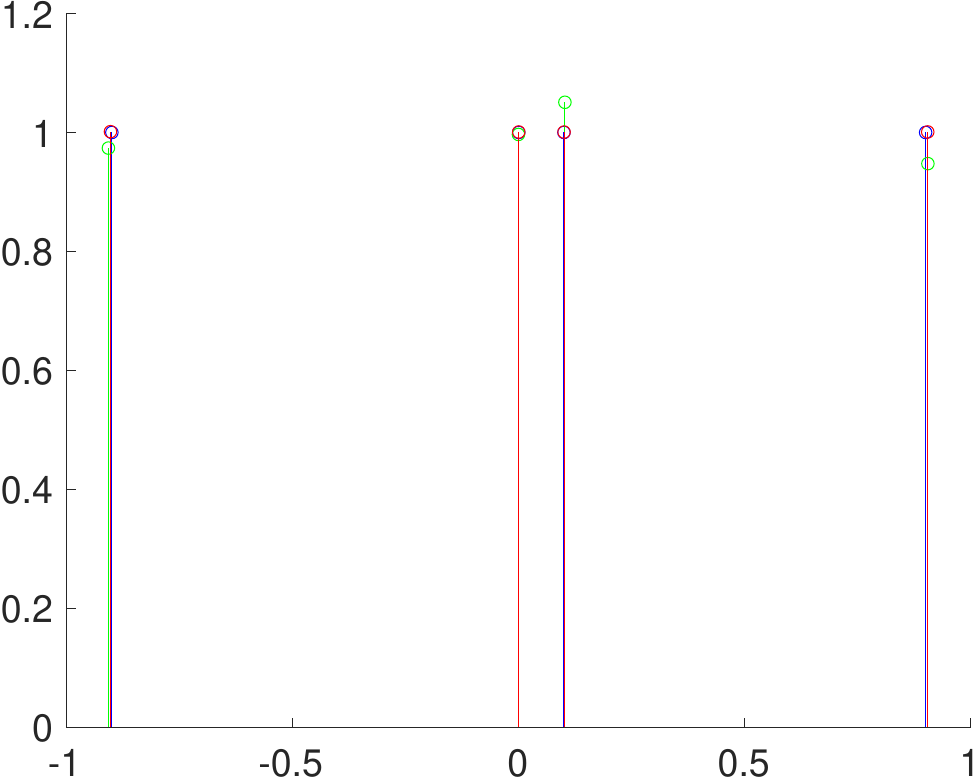}
  \includegraphics[scale=0.25]{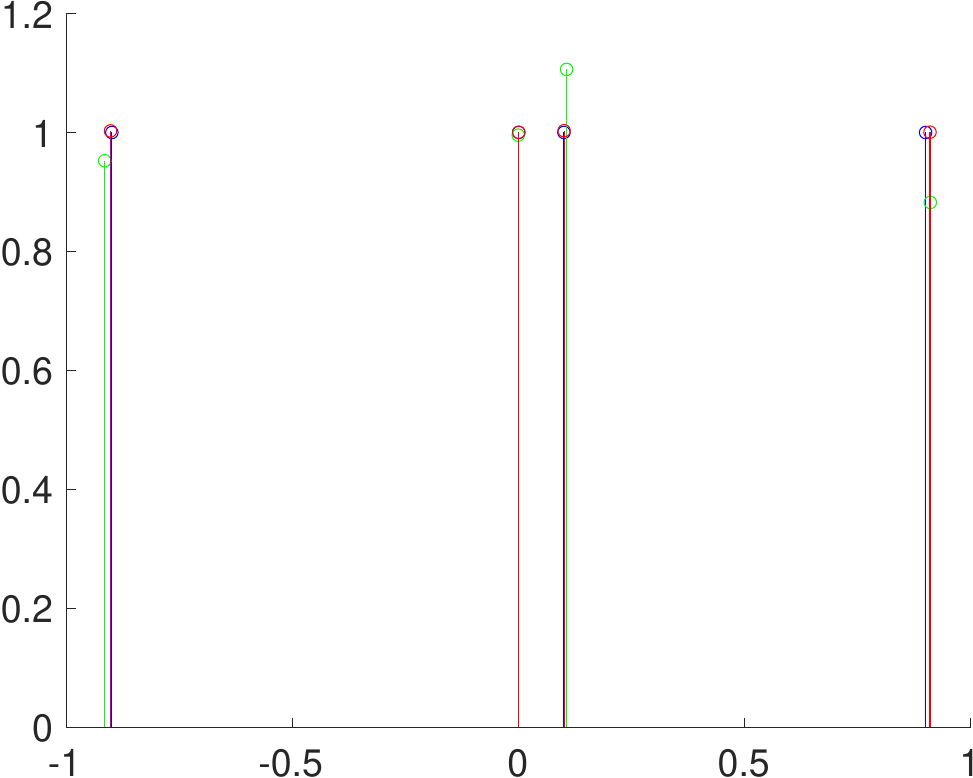} 
  \caption{Spectral function approximation. $G(s,x) = \frac{1}{s-x}$. $X=[-1,1]$. $\{s_j\}$ is the Matsubara grid from $-\frac{(2N-1)\pi}{\beta}i$ to $\frac{(2N-1)\pi}{\beta}i$ with $\beta=100$ and $N=128$. Columns: $\sigma$ equals to $10^{-2}$, $10^{-3}$, and $10^{-4}$. Rows: the easy test (well-separated) and the hard test (with two nearby spikes).}
  \label{fig:exS}
  \end{figure}
  
  Figure \ref{fig:exS} summarizes the experimental results. $n_a=32$. The three columns correspond to $\sigma$ equal to $10^{-2}$, $10^{-3}$, and $10^{-4}$, respectively. Two tests are performed. In the first one (top row), $\{x_k\}$ are well-separated. The reconstructions are accurate for all $\sigma$ values. In the second one (bottom row), two of the spike locations are within $0.1$ distance from each other. In this harder case, the reconstructions also remain accurate for all $\sigma$ values. Notice that the eigenmatrix provides a sufficiently accurate initial guess for postprocessing.
\end{example}

%-------
\begin{example}[Fourier inversion] The problem setup is
  \begin{itemize}
  \item $G(s,x) = \exp(\pi i s x)$.
  \item $X=[-1,1]$.
  \item $\{s_j\}$ are randomly chosen points in $[-5,5]$. $n_s=128$.
  \end{itemize}

  \begin{figure}[h!]
  \centering
  \includegraphics[scale=0.25]{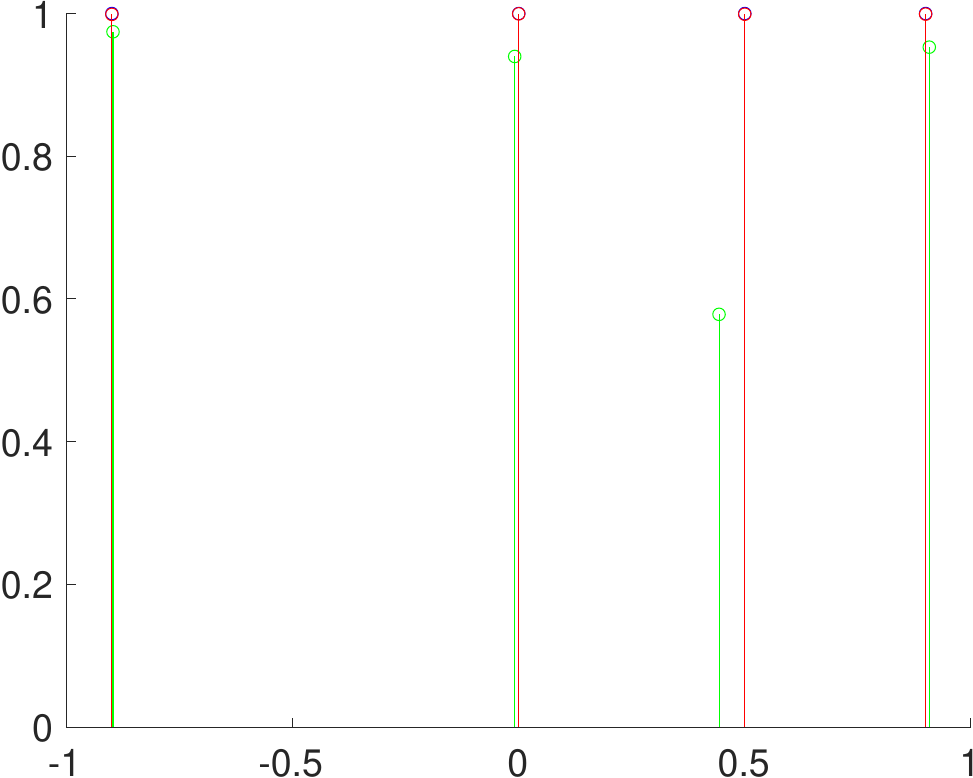}
  \includegraphics[scale=0.25]{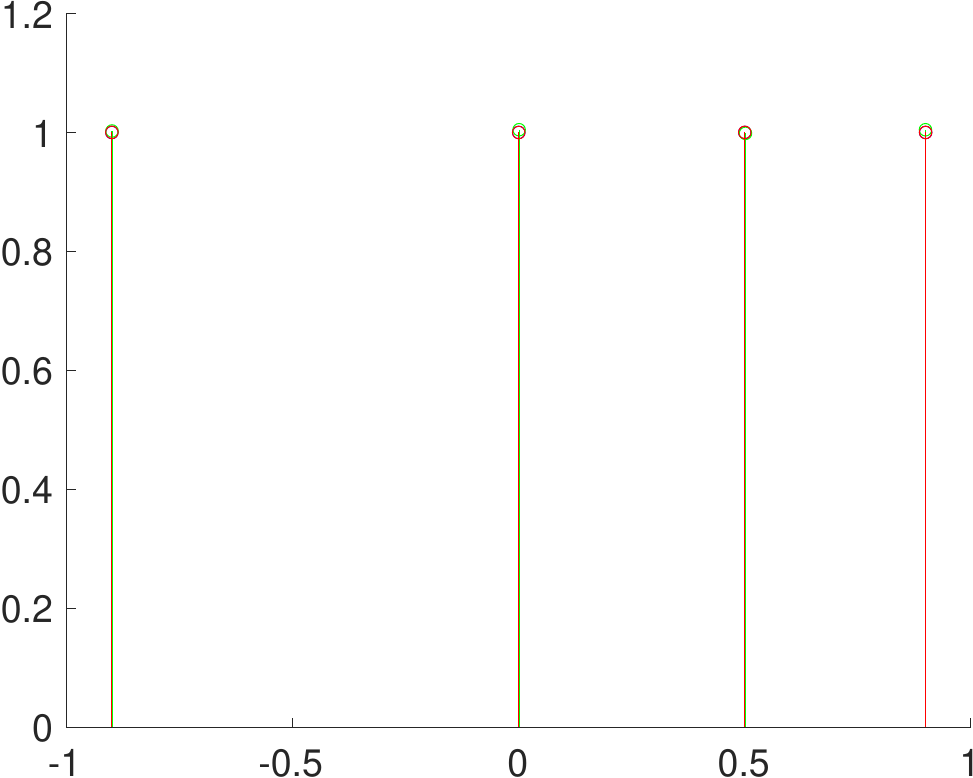}
  \includegraphics[scale=0.25]{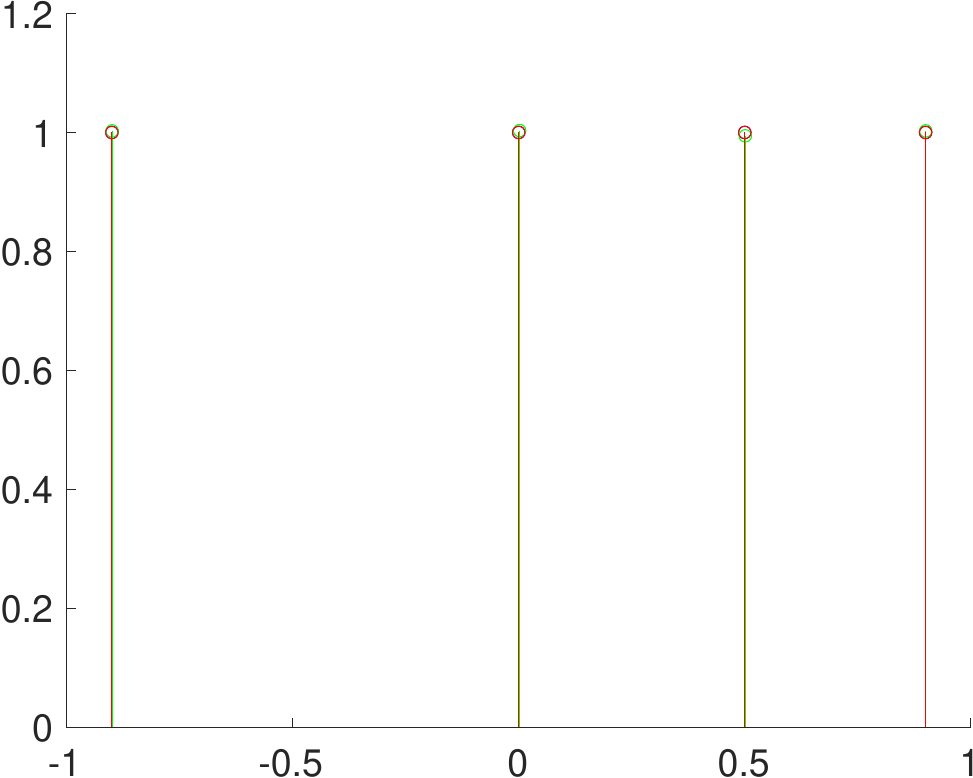}\\
  \includegraphics[scale=0.25]{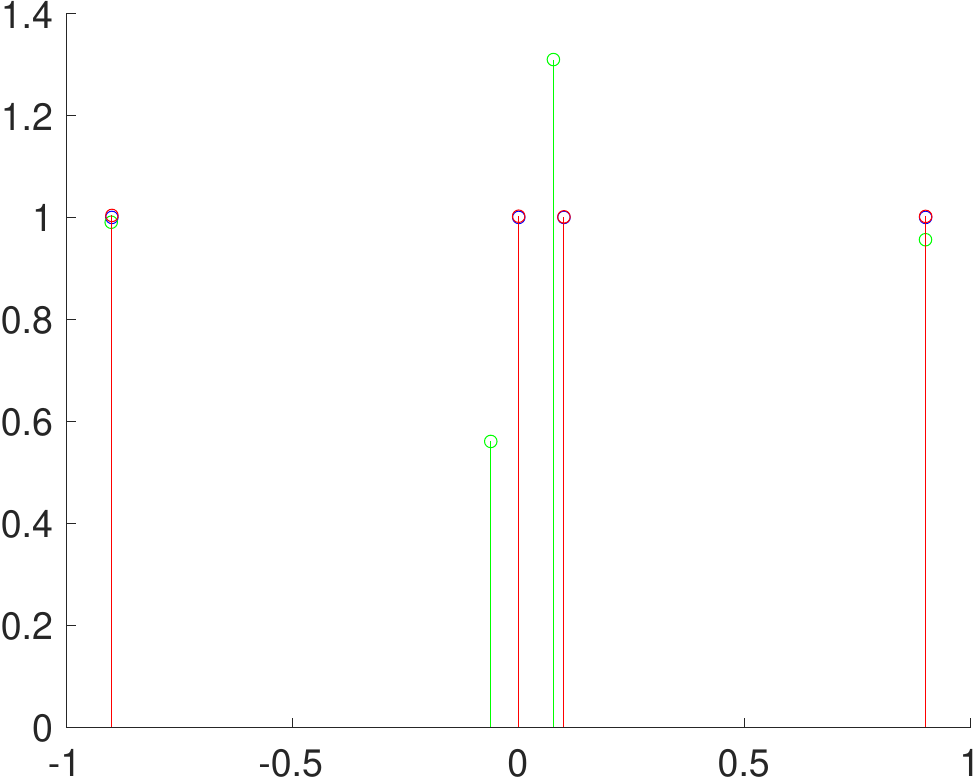}
  \includegraphics[scale=0.25]{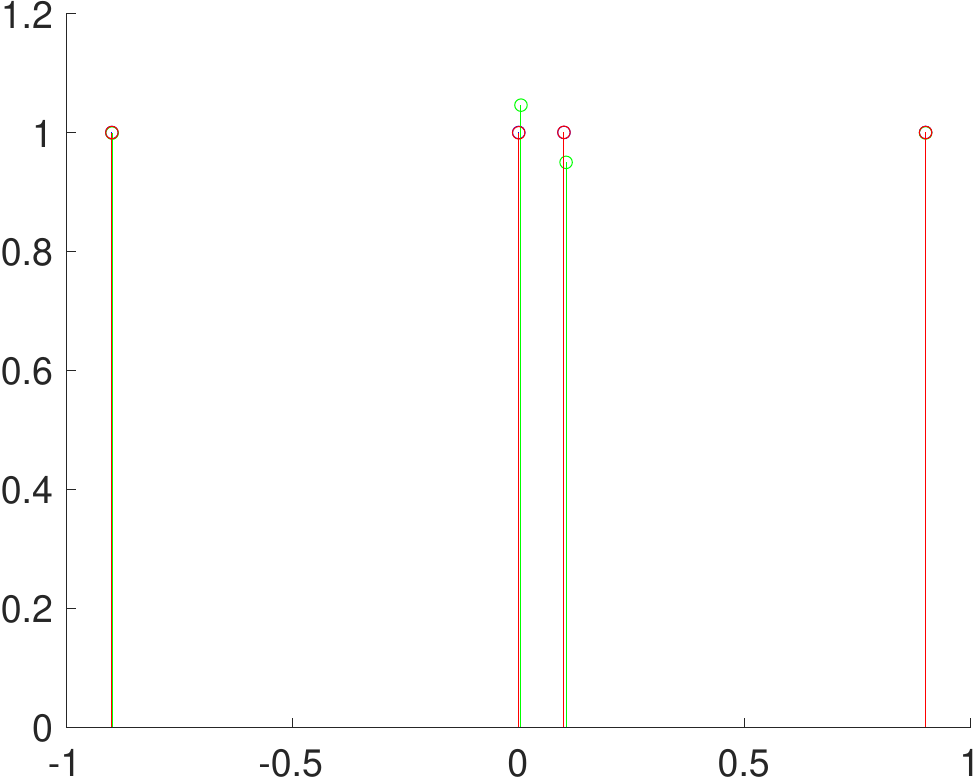}
  \includegraphics[scale=0.25]{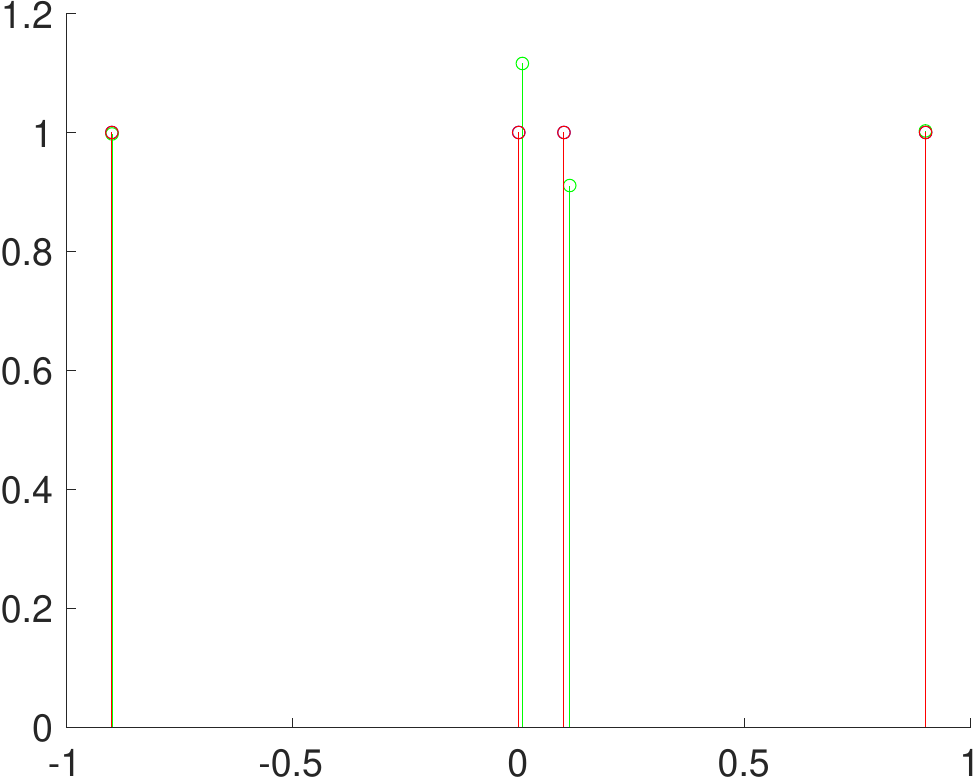} 
  \caption{Fourier inversion. $G(s,x) = \exp(\pi i s x)$. $X=[-1,1]$. $\{s_j\}$ are randomly chosen points in $[-5,5]$. $n_s=128$. Columns: $\sigma$ equals to $10^{-2}$, $10^{-3}$, and $10^{-4}$, respectively. Rows: the easy test (well-separated) and the hard test (with two nearby spikes).}
  \label{fig:exF}
  \end{figure}

  Figure \ref{fig:exF} summarizes the experimental results. $n_a=32$. The three columns correspond to $\sigma$ equal to $10^{-2}$, $10^{-3}$, and $10^{-4}$.  Two tests are performed. In the first one (top row), $\{x_k\}$ are well-separated. The reconstructions are accurate for all $\sigma$ values. In the second one (bottom row), two of the spike locations are within $0.1$ distance from each other. The reconstructions are also accurate for all $\sigma$ values. The eigenmatrix is again able to provide sufficient accurate initial guesses for the postprocessing step.

\end{example}

%-------
\begin{example}[Laplace inversion] The problem setup is
  \begin{itemize}
  \item $G(s,x) = x\exp(-sx)$.
  \item $X=[0.1,2.1]$.
  \item $\{s_j\}$ are random samples in $[0,10]$. $n_s=100$. 
  \end{itemize}
  
  \begin{figure}[h!]
  \centering
  \includegraphics[scale=0.25]{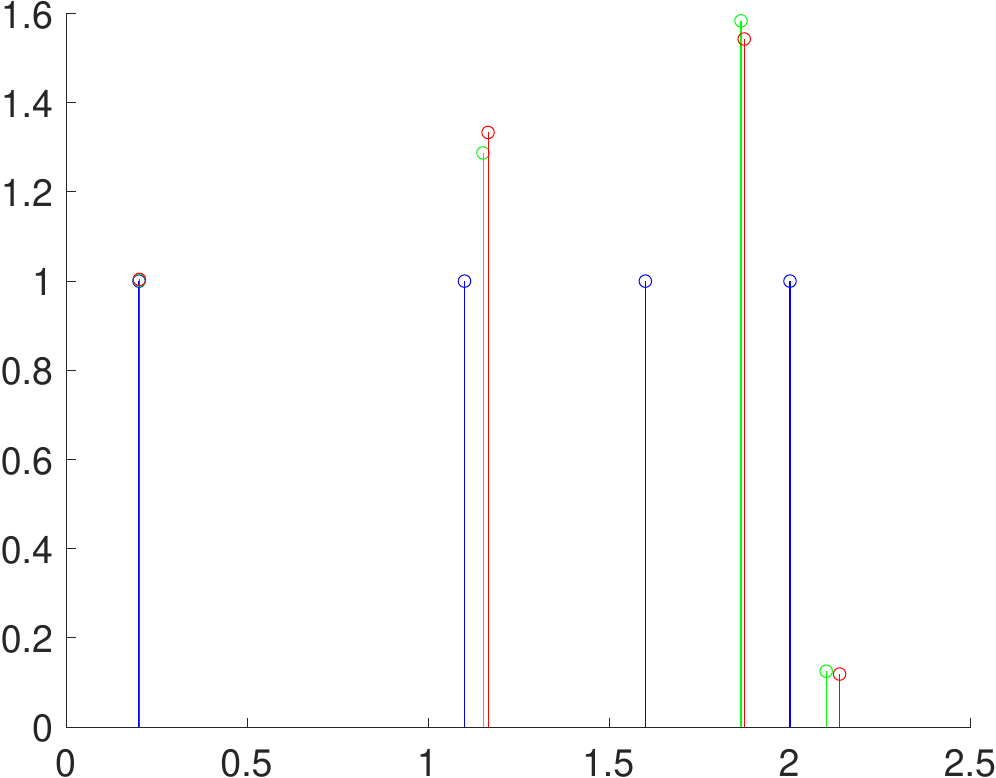}
  \includegraphics[scale=0.25]{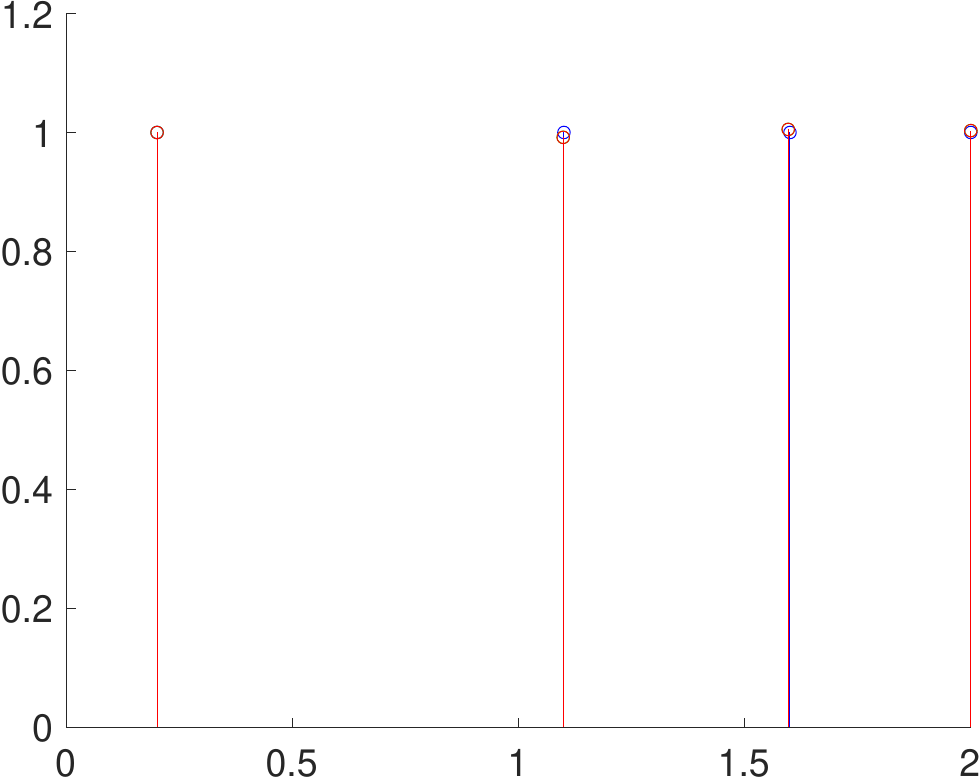}
  \includegraphics[scale=0.25]{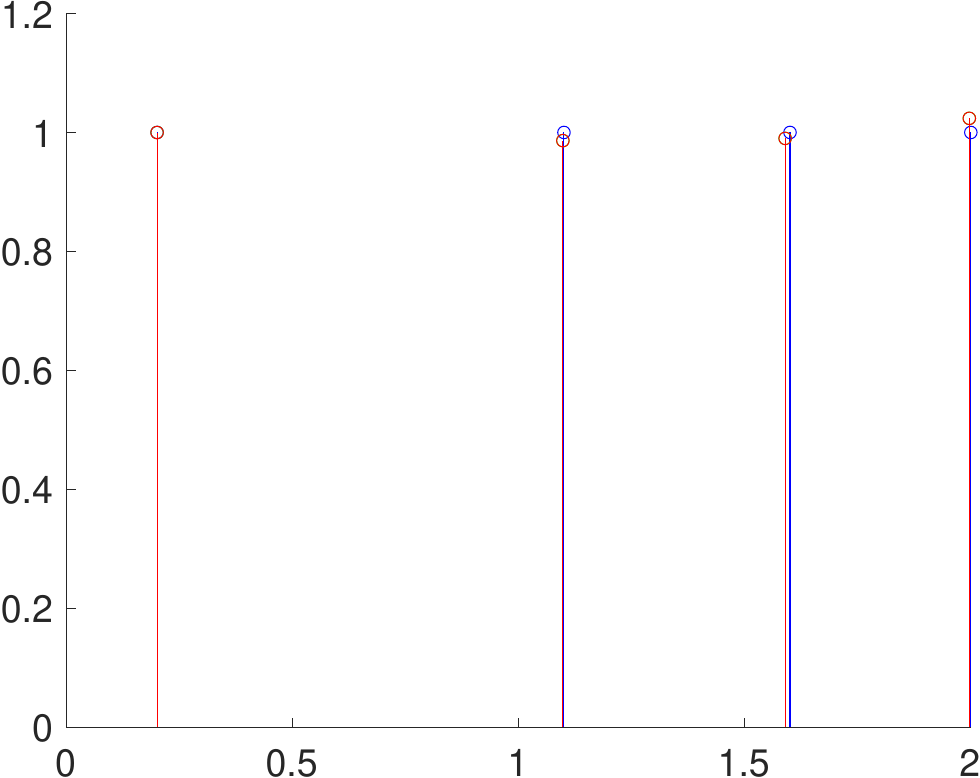}\\
  \includegraphics[scale=0.25]{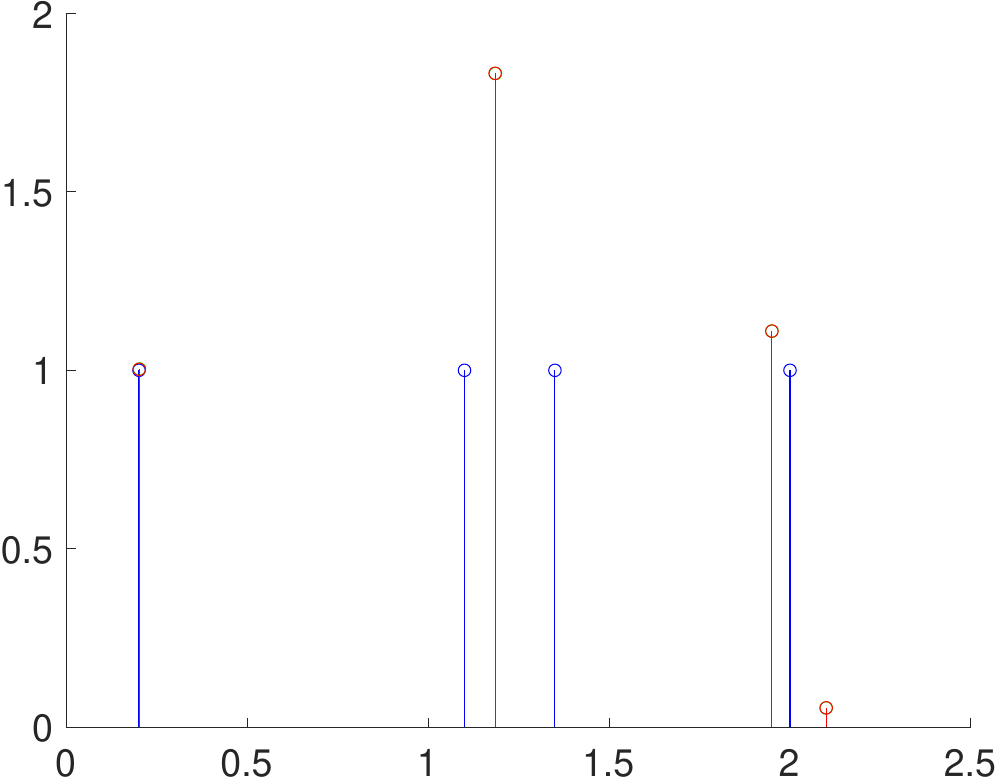}
  \includegraphics[scale=0.25]{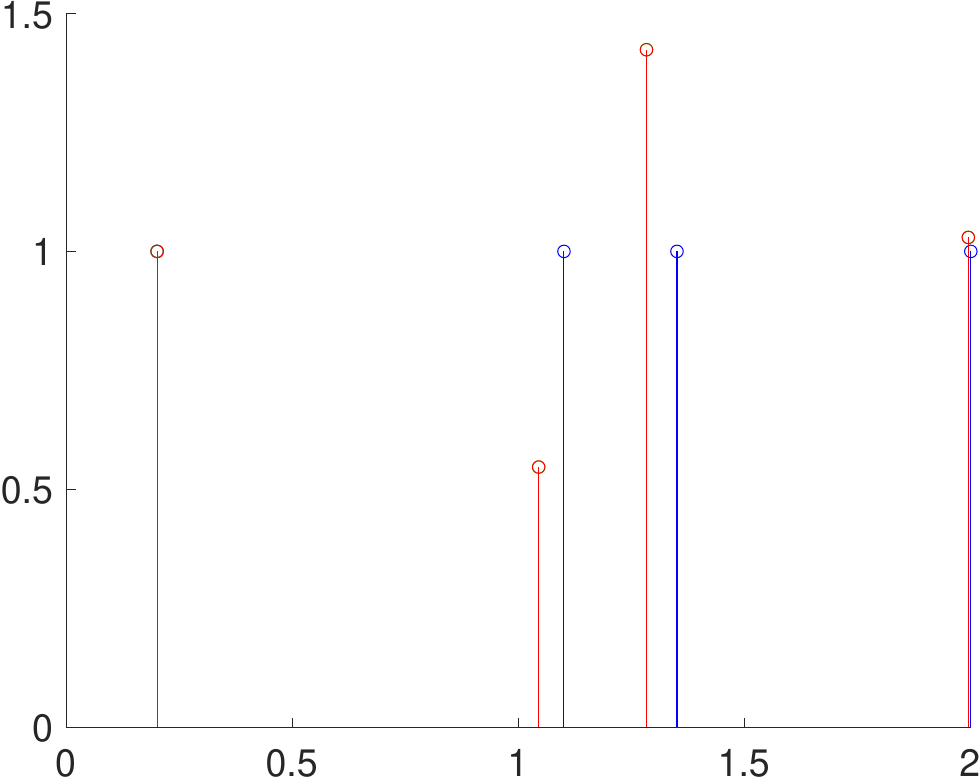}
  \includegraphics[scale=0.25]{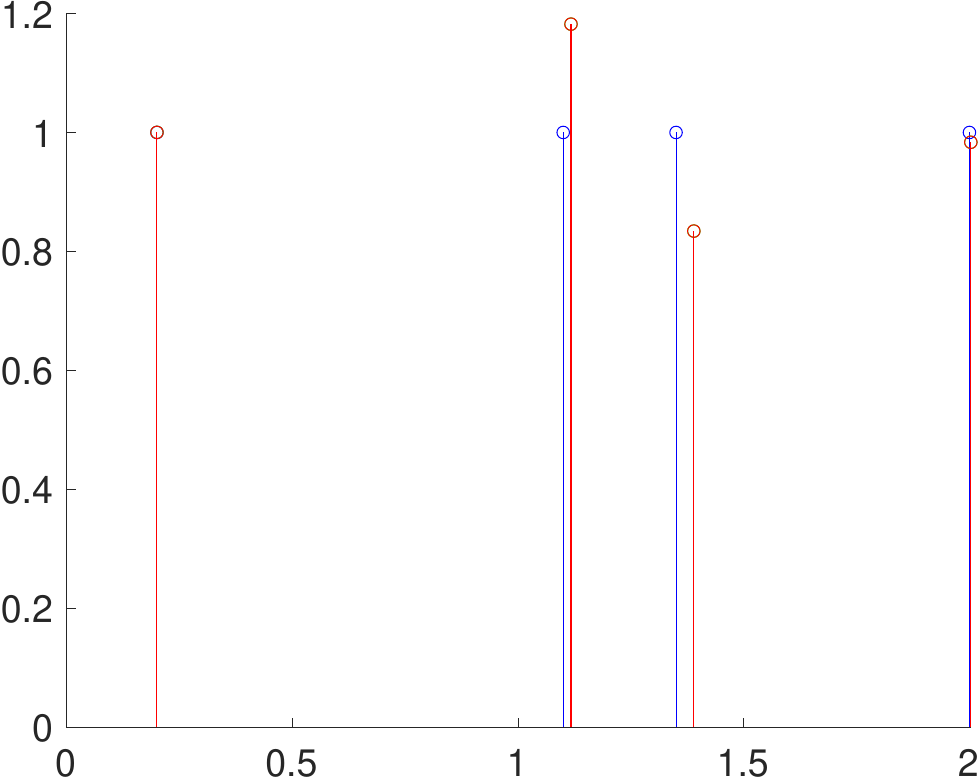} 
  \caption{Laplace inversion.  $G(s,x) = \exp(-sx)$.  $X=[0.1,2.1]$.  $\{s_j\}$ are random samples in $[0,10]$. $n_s=100$. Columns: $\sigma$ equals to $10^{-5}$, $10^{-6}$, and $10^{-7}$, respectively.  Rows: the easy test (well-separated) and the hard test (with two nearby spikes).}
  \label{fig:exL}
  \end{figure}
  
  Figure \ref{fig:exL} summarizes the experimental results. $n_a=32$. The inverse Laplace transform is well-known for its sensitivity to noise. As a result, significantly smaller noise magnitudes are used in this example: the three columns correspond to $\sigma$ equal to $10^{-5}$, $10^{-6}$, and $10^{-7}$. Two tests are performed. In the first one (top row), $\{x_k\}$ are well-separated. The reconstructions are acceptable for $\sigma=10^{-6}$ and accurate for $\sigma=10^{-7}$. In the second harder test (bottom row), two of the spike locations are within $0.25$ distance from each other. The reconstructions provide reasonable reconstructions at $\sigma=10^{-7}$, but significant errors for larger $\sigma$ values.
\end{example}

%-------
\begin{example}[Sparse deconvolution] The problem setup is 
  \begin{itemize}
  \item $G(s,x) = \frac{1}{1+\gamma(s-x)^2}$ with $\gamma=4$.
  \item $X=[-1,1]$.
  \item $\{s_j\}$ are random samples from $[-5,5]$. $n_s=100$.
  \end{itemize}
  
  \begin{figure}[t!]
  \centering
  \includegraphics[scale=0.25]{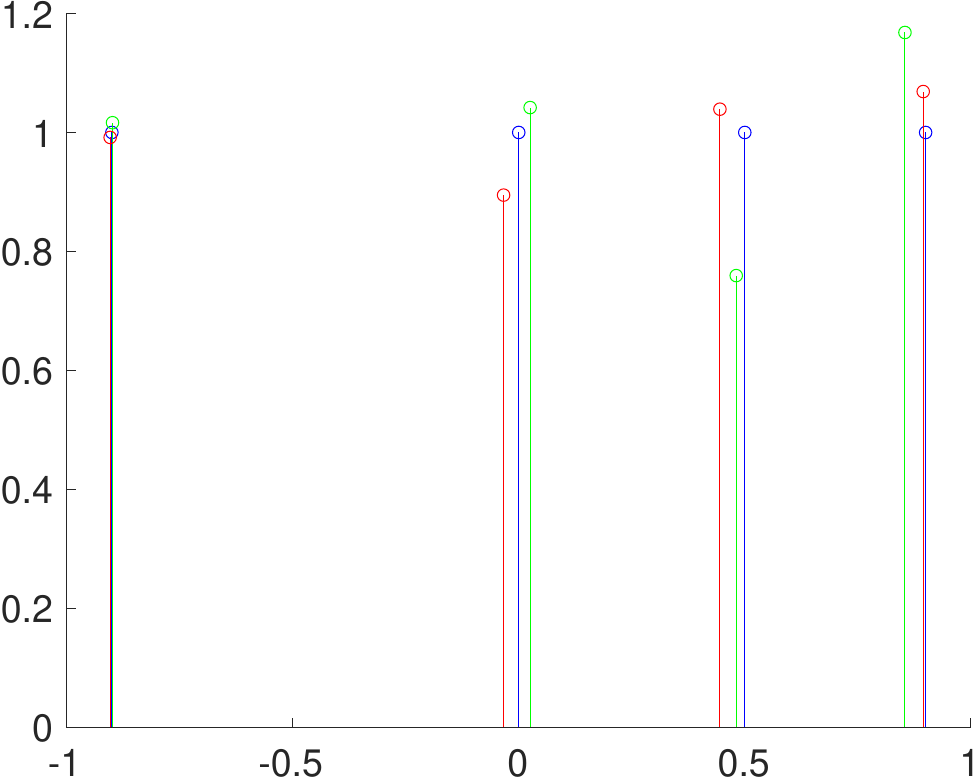}
  \includegraphics[scale=0.25]{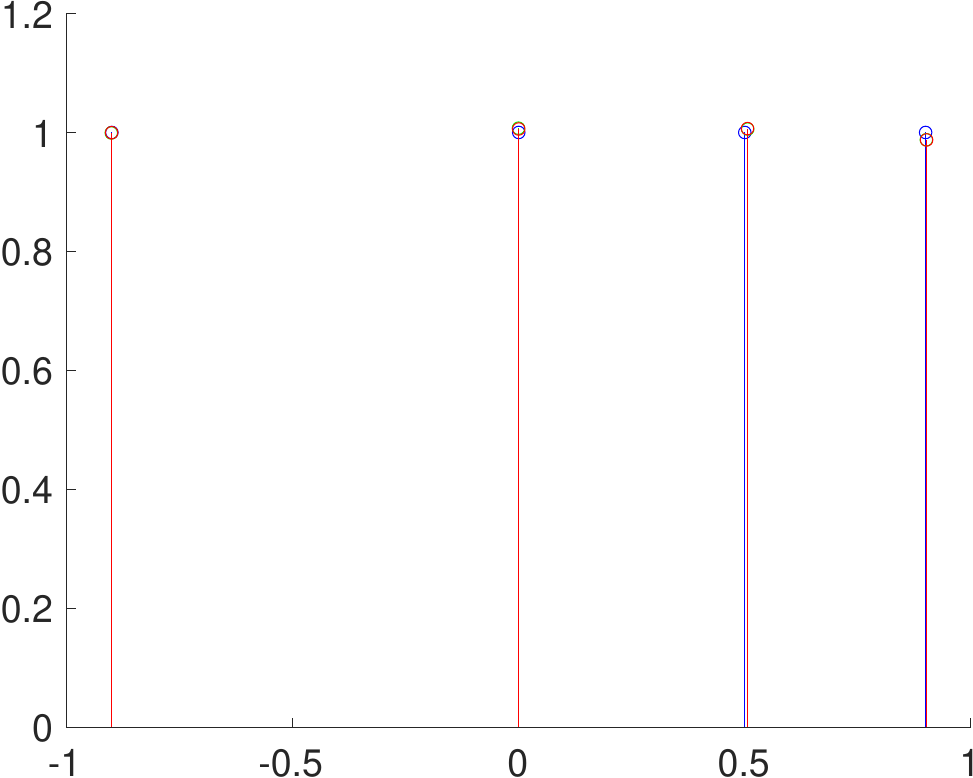}
  \includegraphics[scale=0.25]{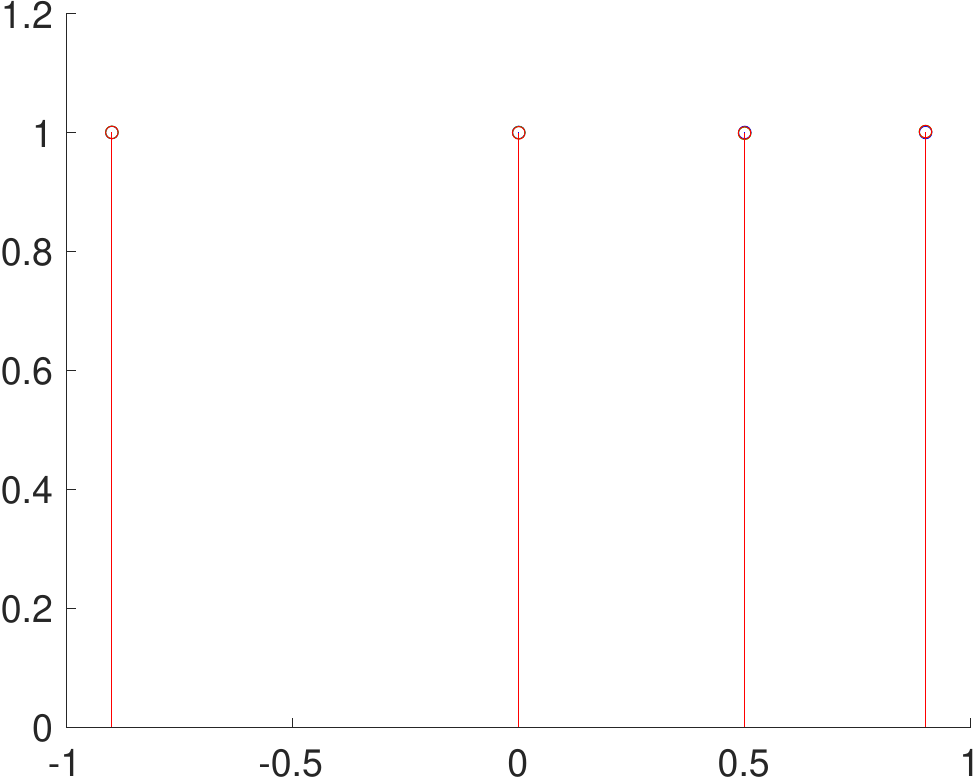}\\
  \includegraphics[scale=0.25]{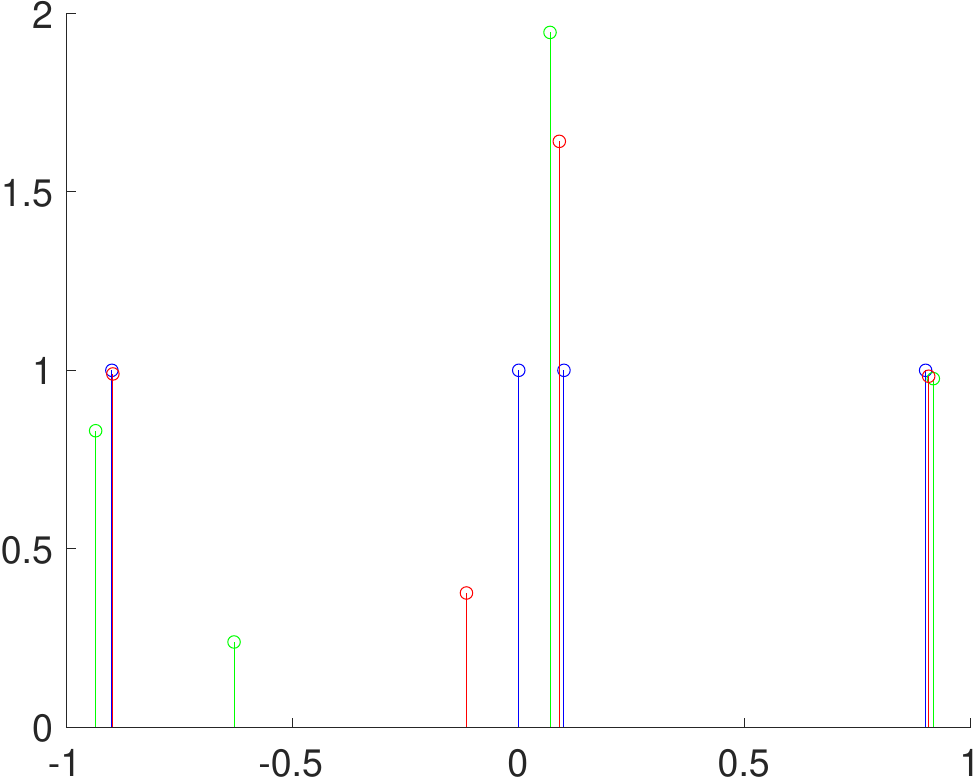}
  \includegraphics[scale=0.25]{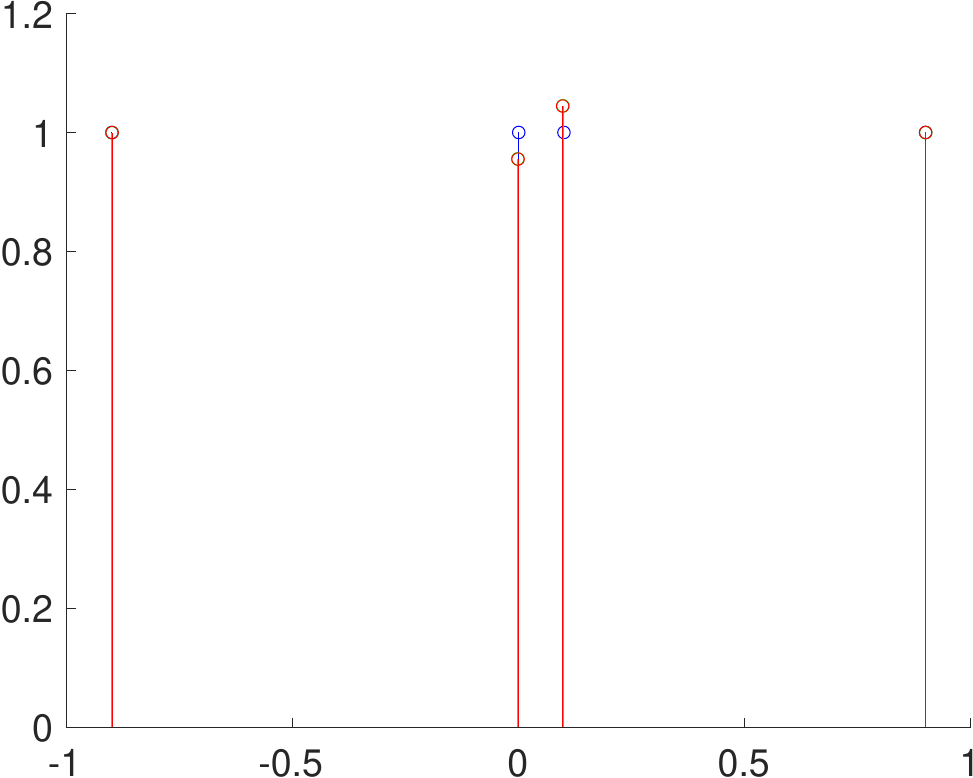}
  \includegraphics[scale=0.25]{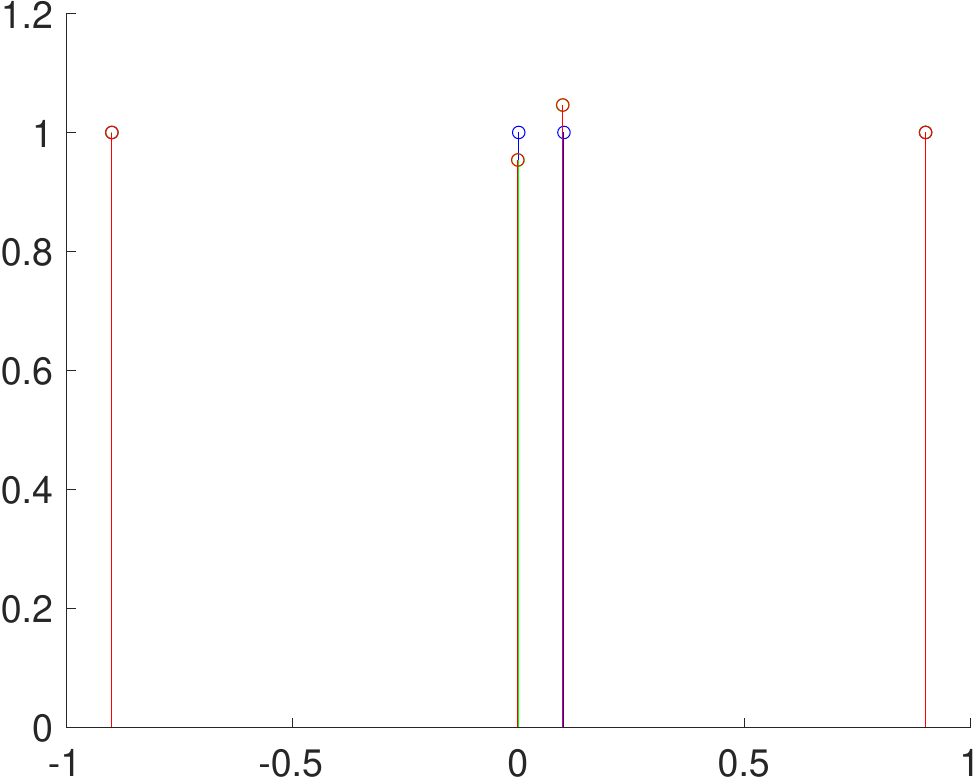} 
  \caption{Sparse deconvolution. $G(s,x) = \frac{1}{1+\gamma(s-x)^2}$ with $\gamma=4$.  $X=[-1,1]$.  $\{s_j\}$ are random samples from $[-5,5]$. $n_s=100$. Columns: $\sigma$ equals to $10^{-2}$, $10^{-3}$, and $10^{-4}$, respectively. Rows: the easy test (well-separated) and the hard test (with two nearby spikes).}
  \label{fig:exD}
  \end{figure}
  
  Figure \ref{fig:exD} summarizes the experimental results. $n_a=32$. The three columns correspond to $\sigma$ equal to $10^{-2}$, $10^{-3}$, and $10^{-4}$. Two tests are performed. In the first one (top row), $\{x_k\}$ are well-separated. The reconstructions are reasonable for $\sigma=10^{-2}$ and accurate for the smaller $\sigma$ values. In the second one (bottom row), two of the spike locations are within $0.1$ distance from each other. The reconstructions are accurate for $\sigma$ equal to $10^{-3}$ and $10^{-4}$.
\end{example}

\begin{remark}
  The numerical experience suggests two lessons important for accurate reconstruction. First, it is important to fully exploit the prior information about the support of the spikes, i.e., making the candidate parameter set $X$ as compact as possible. Second, using the Chebyshev grid (for the real case) and the uniform grid (for the complex case) ensures that $M\bg(x) \approx x \bg(x)$ up to a high accuracy numerically.
\end{remark}

%----------------------------------------------------------
\section{Discussions}\label{sec:disc}

This note introduces the eigenmatrix construction for unstructured sparse recovery problems. It assumes no structure on the sample locations and offers a rather unified framework for such sparse recovery problems. This note is only an exploratory study of the data-driven approach for unstructured sparse recovery, and there are several clear directions for future work.
\begin{itemize}
\item Providing a more principled way of choosing the size of the grid $\{a_t\}$ and the thresholding value for computing $M$.
\item Providing the error estimates of the eigenmatrix approach for the problems mentioned in Section \ref{sec:intro}.
\item Once the eigenmatrix is constructed, the recovery algorithm presented above follows Prony's method and the ESPRIT algorithm. An immediate extension is to combine the eigenmatrix with other algorithms, such as MUSIC and the matrix pencil method.
\end{itemize}

\bibliographystyle{abbrv}

\bibliography{ref}

\end{document}